\documentclass[aos,MSNbibl,seceqn,dvips]{arximspdf}
\usepackage{graphicx}
\doi{10.1214/14-AOS1256} 
\volume{42}
\issue{5}
\pubyear{2014}
\firstpage{2092}
\lastpage{2137}
\docsubty{FLA}

\makeatletter
\newcommand{\SML}{\mathrm{SML}}
\newcommand{\MSLE}{\mathrm{MSLE}}
\newcommand{\IC}{\mathrm{IC}}
\newcommand{\CS}{\mathrm{CS}}
\newcommand{\toy}{\mathrm{toy}}

\newcommand{\rrvert}{\vert}

\newcommand{\rrVert}{\Vert}
\newcommand{\llvert}{\vert}
\newcommand{\llVert}{\Vert}
\def\IK{\mathbb{K}}
\newtheorem{theorem}{Theorem}[section]
\newtheorem{corollary}{Corollary}[section]
\newtheorem{lemma}{Lemma}[section]
\newproclaim{remark}{Remark}[section]
\newtheorem{condition}{Condition}[section]
\newcommand{\fraca}[2]{{#1}/{#2}}

\newcommand{\fracc}[2]{{#1}/(#2)}

\makeatother

\begin{document}
\begin{frontmatter}

\title{Maximum smoothed likelihood estimators\break for the interval
censoring model}
\runtitle{The MSLE for interval censoring}

\begin{aug}
\author{\fnms{Piet}~\snm{Groeneboom}\corref{}\ead[label=e1]{P.Groeneboom@tudelft.nl}\ead[label=u1,url]{http://dutiosc.twi.tudelft.nl/\textasciitilde pietg/}}
\runauthor{P. Groeneboom}
\affiliation{Delft University}
\address{Delft Institute of Applied Mathematics\\
Delft University of Technology\\
Mekelweg 4, 2628 CD Delft\\
The Netherlands\\
\printead{e1}\\
\printead{u1}}
\end{aug}

\received{\smonth{3} \syear{2012}}
\revised{\smonth{7} \syear{2014}}

%
\begin{abstract}
We study the maximum smoothed likelihood estimator (MSLE) for interval
censoring, case 2, in the so-called separated case. Characterizations
in terms of convex duality conditions are given and strong consistency
is proved. Moreover, we show that, under smoothness conditions on the
underlying distributions and using the usual bandwidth choice in
density estimation, the local convergence rate is $n^{-2/5}$ and the
limit distribution is normal, in contrast with the rate $n^{-1/3}$ of
the ordinary maximum likelihood estimator.
\end{abstract}

%
\end{frontmatter}

\section{Introduction}\label{sectionintro}
In \cite{pietgeurtbirgit10}, the maximum smoothed likelihood
estimator (MSLE) and smoothed maximum likelihood estimator (SMLE) were
studied for the current status model, the simplest interval censoring
model. It is called the interval censoring, case 1, model in \cite
{Gr91} and \cite{GrWe92}. It was shown in \cite
{pietgeurtbirgit10} that, under certain regularity conditions, the
MSLE and the SMLE, evaluated at a fixed interior point, converge at
rate $n^{-2/5}$ to the real underlying distribution function, if one
takes a bandwidth of order $n^{-1/5}$. This convergence rate is faster
than the convergence rate of the nonsmoothed maximum likelihood
estimator, which is $n^{-1/3}$ in this situation, as shown in \cite
{Gr91} and \cite{GrWe92}. Moreover, the limit distribution is
normal, in contrast with the limit distribution of the nonsmoothed
maximum likelihood estimator.

The interval censoring model, where there is an interval in which the
relevant (unobservable) event takes place, is more common, in
particular in medical statistics. It is called the interval censoring,
case 2, model in \cite{Gr91} and \cite{GrWe92}. A preliminary
discussion of the SMLE in this situation can be found in \cite
{piettom11}, where it was shown that the development of the theory of
the SMLE for this model crucially depends on a further analysis of the
integral equations, studied in \cite{GeGr99,GeGr96} and \cite{GeGr97}. In the present paper, we study the
MSLE and prove a consistency and asymptotic normality result for this
estimator. We also discuss algorithms for computing the MSLE, which is
a rather complicated issue.

We recall the interval censoring, case 2, model. Let $X_1,\ldots,X_n$
be a sample of unobservable random variables
from an unknown distribution function $F_0$ on $[0,\infty)$.
Suppose that one can observe $n$ pairs $(T_i,U_i)$, independent of
$X_i$, where $U_i>T_i$. Moreover,
\[
\Delta_{i1}\stackrel{\mathrm{def}}=1_{\{X_i\le T_i\}},\qquad \Delta
_{i2}\stackrel{\mathrm{def}}=1_{\{T_i<X_i\le U_i\}}\quad\mbox{and}\quad\Delta
_{i3}\stackrel{\mathrm{def}}=1-\Delta_{i1}-\Delta
_{i2}, %
\]
provide the only information one has on the position of the random variables
$X_i$ with respect to the observation times $T_i$ and $U_i$. In
this set-up, one wants to estimate the unknown distribution function
$F_0$, generating the ``unobservables'' $X_i$, on an interval $[0,M]$.

Interestingly, from a computational point of view, the MLE for the
distribution function of the hidden variable in the case that one has
more observation times $T_i, U_i, V_i,\ldots$ ``per hidden variable,''
can always be reduced to the case of interval censoring, case 2. This
follows from the fact that at most two of the observation times of the
set $\{T_i, U_i, V_i,\ldots\}$ are relevant for the location of the
hidden variable. If we know that the hidden variable is located between
two observation times, while the other observation times for this
hidden variable are either more to the right or more to the left, then
these other observation times do not give extra information and can be
discarded in computing the MLE. Likewise, if we know that the hidden
variable lies to the right of all these observation times, all
observation times smaller than the largest one do not give extra
information, with a similar situation if we know that the hidden
observation time lies to the left of the smallest observation time for
this variable. So, in the last two cases, only one observation time
gives relevant information and the other ones can be discarded. This
motivates concentrating on the interval censoring, case 2, model, as an
extension of the current status model.

The MSLE (maximum smoothed likelihood estimator) is defined in the
following way. Let $g$ be the joint density of the observation pairs
$(T_i,U_i)$, with first marginal $g_1$ and second marginal $g_2$.
Moreover, let the densities $h_{01}$, $h_{02}$ and $h_0$ be defined by
%
\begin{eqnarray}\label{defh0}
h_{01}(t)&=&F_0(t)g_1(t),\nonumber
\\
h_{02}(u)&=&\bigl\{1-F_0(u)\bigr\}g_2(u),
\\
h_0(t,u)&=&\bigl\{F_0(u)-F_0(t)\bigr\}g(t,u).\nonumber
\end{eqnarray}
We define $\tilde h_{nj}$, $j=1,2$ and $\tilde h_n$ as the estimates of
the densities $h_{0j}$, $j=1,2$ and the 2-dimensional density $h_0$,
respectively, where
%
\begin{eqnarray}
\label{defhnj} \tilde h_{n1}(t)&=&\frac{1}n\sum
_{i=1}^nK_{b_n}(t-T_i) \Delta
_{i1},\qquad\tilde h_{n2}(u)=\frac{1}n\sum
_{i=1}^nK_{b_n}(u-U_i)\Delta
_{i3},
\\
%
%
\label{defhn} \qquad\tilde h_n(t,u)&=&\frac{1}n\sum
_{i=1}^nK_{b_n}(t-T_i)
K_{b_n}(u-U_i)\Delta_{i2}
\end{eqnarray}
and
\[
K_{b_n}(x)=\frac{1}{b_n}K \biggl(\frac{x}{b_n} \biggr),
\]
for a symmetric continuously differentiable kernel $K$ with compact
support, like the triweight kernel
%
\begin{equation}
\label{defK} K(x)=\tfrac{35}{32} \bigl(1-x^2
\bigr)^31_{[-1,1]}(x),\qquad x\in\mathbb R.
\end{equation}
At boundary points, we use a boundary correction by replacing the
kernel $K$ by a linear combination of $K(u)$ and $uK(u)$. For example,
if $t\in[0,b_n)$, we define
\[
\tilde h_{n1}(t)=\alpha(t/b_n)\frac{1}n\sum
_{i=1}^nK_{b_n}(t-T_i)
\Delta_{i1} +\beta(t/b_n)\frac{1}n\sum
_{i=1}^n \frac{t-T_i}{b_n}K_{b_n}(t-T_i)
\Delta_{i1}, %
\]
where the coefficients $\alpha(u)$ and $\beta(u)$ are defined by
%
\begin{equation}
\label{coeff1} \alpha(u)\int_{-1}^u K(x) \,dx+
\beta(u)\int_{-1}^u xK(x) \,dx=1,\qquad u\in[0,1]
\end{equation}
and
%
\begin{equation}
\label{coeff2} \alpha(u)\int_{-1}^u xK(x) \,dx+
\beta(u)\int_{-1}^u x^2K(x) \,dx=0,\qquad u
\in[0,1].
\end{equation}
It may happen that $\tilde h_{n1}(t)<0$; in that case we put $\tilde
h_{n1}(t)=0$.

If $t\in(M-b_n,M]$, we similarly define
\begin{eqnarray*}
\tilde h_{n1}(t)&=&\alpha\bigl((M-t)/b_n\bigr)
\frac{1}n\sum_{i=1}^nK_{b_n}(t-T_i)
\Delta_{i1}
\\
&&{} -\beta\bigl((M-t)/b_n\bigr)\frac{1}n
\sum_{i=1}^n \frac
{t-T_i}{b_n}K_{b_n}(t-T_i)
\Delta_{i1}, %
\end{eqnarray*}
where the functions $\alpha$ and $\beta$ are again defined by (\ref
{coeff1}) and (\ref{coeff2}).
The estimates $\tilde h_{n2}$ and $\tilde h_n$ are similarly defined if
one or more (in the case of $\tilde h_n$) arguments have distance less
than $b_n$ to the boundary; for $\tilde h_n$ we apply this to the
factors of the product of the kernels separately, in the same way as
for the one-dimensional estimates $\tilde h_{nj}$. We finally divide
$\tilde h_{n1}(t)$, $\tilde h_{n2}(t)$ and $\tilde h_n(t,u)$ by
\[
\int_{[0,M]} \bigl\{\tilde h_{n1}(x)+\tilde
h_{n2}(x) \bigr\} \,dx +\int_{[0,M]^2}\tilde
h_n(x,y) \,dx \,dy, %
\]
(i.e., by a discrete approximation to this quantity) to give a total
mass approximately equal to $1$ to the observation density.

The MSLE $\hat F_n$ is now defined as the distribution function,
maximizing the criterion function
%
\begin{eqnarray}\label{criterionfunction1}
\ell(F) &=& \int\tilde h_{n1}(t)\log F(t) \,dt+
\int\tilde
h_{n2}(u)\bigl\{1-F(t)\bigr\} \,du
\nonumber\\[-8pt]\\[-8pt]\nonumber
&&{} +\int\tilde h_n(t,u)\log
\bigl\{ F(u)-F(t)\bigr\} \,dt \,du,
\end{eqnarray}
as a function of the distribution function $F$. But in practice we
discretize, and maximize
%
\begin{eqnarray}
\label{criterionfunction2} &&\sum_{i=1}^m \bigl\{\tilde
h_{n1}(t_i)\log F(t_i) \bigr
\}d_i +\sum_{i=1}^m \bigl\{
\tilde h_{n2}(t_i)\log\bigl\{1-F(t_i)\bigr\}
\bigr\} d_i
\nonumber\\[-8pt]\\[-8pt]\nonumber
&&\qquad{}+\sum_{i=1}^{m-1}
\sum_{j=i+1}^m \bigl\{ \tilde
h_n(t_i,t_j)\log\bigl\{F(t_j)-F(t_i)
\bigr\} \bigr\}d_id_j,
\end{eqnarray}
over all distribution functions $F$, where $0=t_0<t_1,\ldots,<t_m=M$ are
the points of a grid and $d_i=t_i-t_{i-1}$, $i=1,\ldots,m$.

Note that $\ell(F)$ is a smoothed log likelihood for $F$ and,
therefore, the maximizing (sub)distribution function $F$ is called the
maximum smoothed likelihood estimator (MSLE).
Also note that the maximization of (\ref{criterionfunction1}) is the
same as the minimization of the Kullback--Leibler distance
\begin{eqnarray*}
&& \int\tilde h_{n1}(t)\log\frac{\tilde h_{n1}(t)}{F(t)\tilde
g_{n1}(t)} \,dt +\int\tilde
h_{n1}(t)\log\frac{\tilde h_{n2}(t)}{\{1-F(t)\}\tilde
g_{n2}(t)} \,dt
\\
&&\qquad{}+\int\tilde h_n(t,u)\log\frac{\tilde h_n(t,u)}{\{F(u)-F(t)\}
\tilde g_n(t,u)} \,dt \,du,
\end{eqnarray*}
where $\tilde g_{ni}$ and $\tilde g_n$ are kernel estimates of the
densities $g_i$ and $g$, computed in the same way as the estimates
$h_{ni}$ and $\tilde h_n$ (but without the indicators $\Delta_{ij}$).

\begin{figure}

\includegraphics{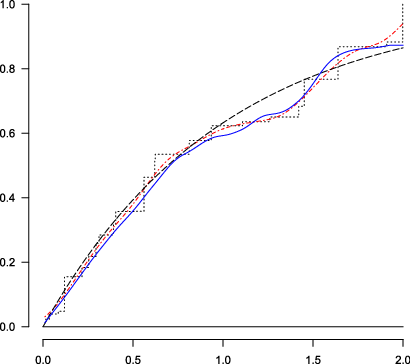}

\caption{The MSLE (solid), SMLE (dashed-dotted) and MLE (dotted) on
$[0,1]$ for a sample of size $n=1000$ from the exponential distribution
function $F_0(x)=1-\exp\{- x\}$ (dashed); the bivariate observation
density is $g(x,y)=6(y-x-\varepsilon)^2/\{(2-x-\varepsilon
)(2-\varepsilon)\}^2, x+\varepsilon<y$ on the triangle with
vertices $(0,\varepsilon)$, $(0,2)$ and $(2-\varepsilon,2)$, where
$\varepsilon=0.1$. The bandwidth for the computation of the MSLE\vspace*{1pt} was
$b_n=n^{-1/5}\approx0.25119$.}\label{figSMLEMSLE}
\end{figure}

\begin{figure}

\includegraphics{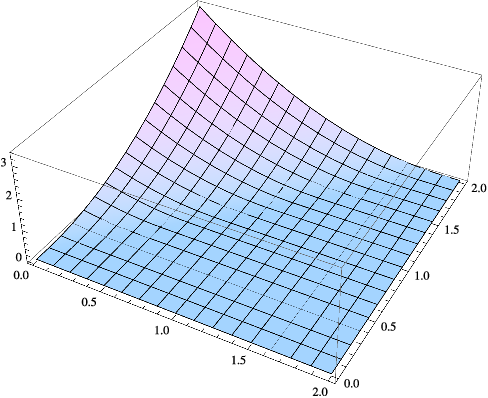}

\caption{The bivariate observation density $g$ on $[0,2]^2$, where
$\varepsilon=0.1$.}
\label{figg}
\end{figure}

Defining the SMLE (smoothed maximum likelihood estimator) is somewhat
easier. If we have computed the ordinary MLE $\hat F_n$, we simply
define the SMLE $F_n^{\SML}(x)$ by
%
\begin{eqnarray}\label{SMLE}
F_n^{\SML}(x)&=&\int\IK_{b_n}(x-y) \,d
\hat F_n(y),
\nonumber\\[-8pt]\\[-8pt]\nonumber
\IK_{b_n}(u)&=&\int_{-\infty}^{u/b_n}K(x)\,dx,
\end{eqnarray}
where $K$ is of type (\ref{defK}) again. A picture of the MLE, the
MSLE and the SMLE for a sample of size $n=1000$ from an exponential
distribution function is shown in Figure~\ref{figSMLEMSLE}. A
picture of the bivariate observation density $g$, with $\varepsilon
=0.1$, is shown in Figure~\ref{figg}.

%

%

\subsection{The SMLE and MSLE for the current status model}\label{subsectioncurstatmodel}

Before embarking on the theory for this model, it might be instructive
to recapitulate the rather different ways in which the asymptotic
distributions of the SMLE and the MSLE are derived for the simpler
current status model. In this case, the data are given by
\[
(T_1,\Delta_1),\ldots,(T_n,\Delta
_n), %
\]
where
\[
\Delta_i=1_{\{X_i\le T_i\}}, %
\]
and $X_i$ and $T_i$ are independent.

Let $\tilde F_n^{(\SML)}$ be the SMLE for the current status model,
defined by (\ref{SMLE}), but now using the MLE $\hat F_n$ in the
current status model. It is shown in \cite{pietgeurtbirgit10} that,
under suitable smoothness conditions, we can write, if $b_n\asymp n^{-1/5}$,
%
\begin{eqnarray}\label{SMLCS}
&& \tilde F_n^{(\SML)}(t)-\int K_{b_n}(t-u) \,dF_0(u)
\nonumber\\[-8pt]\\[-8pt]\nonumber
&&\qquad =\int\theta_{t,b_n,F}^{\CS}(u,
\delta)\,d (\mathbb Q_n-Q_0 ) (u,\delta
)+o_p \bigl(n^{-2/5} \bigr),
\end{eqnarray}
where
%
\begin{equation}
\label{defthetaCS} \theta_{t,b,F}^{\CS}(u,\delta)=-
\frac{\delta\phi
^{\CS}_{t,b,F}(u)}{F(u)}+\frac{(1-\delta)\phi
^{\CS}_{t,b,F}(u)}{1-F(u)},\qquad u\in(0,1)
\end{equation}
and $\phi^{\CS}_{t,b,F}$ is given by
\[
\phi^{\CS}_{t,b,F}(u)=\frac{F(u)\{1-F(u)\}}{g(u)} b^{-1}K
\bigl((t-u)/b\bigr). %
\]
Moreover, $g$ is the density of the (one-dimensional) observation distribution.

The solution $\phi_{t,b_n,F_0}^{\CS}$ gives as an approximation for
$n \operatorname{var}(\tilde F_n(t))$:
\begin{eqnarray*}
{\mathbb E} \theta^{\CS}_{t,b_n,F_0}(T_1,\Delta
_1)^2&=&\int\frac
{\phi_{t,b_n,F_0}^{\CS}(u)^2}{F_0(u)}g(u) \,du +\int
\frac{\phi^{\CS}_{t,b_n,F_0}(u)^2}{1-F_0(u)}g(u) \,du
\\
&\sim&\frac{F_0(t)\{1-F_0(t)\}}{b_ng(t)}\int K(u)^2 \,du,\qquad b_n\to0.
\end{eqnarray*}
Taking the bias into account, we get, if $b_n\asymp n^{-1/5}$, for the
SMLE the central limit theorem
%
\begin{eqnarray}\label{CLCSSMLE}
\qquad \sqrt{n} \biggl\{\tilde F_n^{\CS}(t)-F_0(t)-
\frac12b_n^2f_0'(t)\int
u^2K(u) \,du \biggr\}\Bigm/\sigma_n\stackrel{{\mathcal D}}
\longrightarrow N (0,1 ),
\nonumber\\[-8pt]\\[-12pt]
\eqntext{n\to\infty,}
\end{eqnarray}
where
\[
\sigma_n^2={\mathbb E} \theta^{\CS}_{t,b_n,F_0}(T_1,
\Delta_1)^2\sim\frac{F_0(t)\{1-F_0(t)\}}{b_ng(t)}\int K(u)^2
\,du,\qquad n\to\infty; %
\]
see Theorem~4.2, page 365 \cite{pietgeurtbirgit10}.

\begin{figure}

\includegraphics{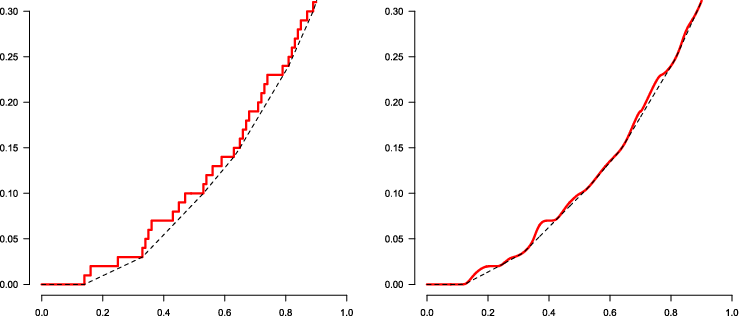}

\caption{Unsmoothed and smoothed cusum diagram.}\label{smoothedunsmoothedcusum}
\end{figure}

On the other hand, for the MSLE in the current status model it is first
shown that the MSLE corresponds to the slope of greatest convex
minorant of the \textit{continuous} cusum diagram
%
\begin{eqnarray}\label{continuouscusum}
&\displaystyle \biggl(\int\IK_b(t-x) \,d{\mathbb G}_n(x),
\int\delta\IK_b(t-x) \,d{\mathbb P}_n(x,\delta) \biggr),&
\nonumber\\[-8pt]\\[-8pt]\nonumber
&\displaystyle \IK_b(y)=\int_{-\infty}^{y/b} K(u) \,du,\qquad t\ge0,&
\end{eqnarray}
where ${\mathbb G}_n$ is the empirical distribution function of the
$T_i$ and ${\mathbb P}_n$ the empirical distribution function of the
pairs $(T_i,\Delta_i)$, analogously to the way the MLE corresponds to
the slope of greatest convex minorant of the cusum diagram
\[
\biggl(\int_{[0,t]}\,d{\mathbb G}_n(u),\int
_{[0,t]}\delta \,d{\mathbb P}_n(u,\delta) \biggr),\qquad t
\ge0. %
\]
A picture of the cusum diagram for the MLE and the SMLE for the same
sample is shown in Figure~\ref{smoothedunsmoothedcusum}.

%

Next, it is shown that the MSLE is at interior points asymptotically
equivalent to the ratio of kernel estimators
%
\begin{equation}
\label{pluginCS} \frac{g_{n,b_n}^{\delta}(t)}{g_{n,b_n}(t)},
\end{equation}
where
\[
g_{n,b_n}(t)=\int K_{b_n}(t-u) \,d{\mathbb G}_n(u),
\qquad g_{n,b_n}^{\delta}(t)=\int\delta K_{b_n}(t-u) \,d{
\mathbb P}_n(u,\delta). %
\]
This leads to the following central limit theorem for the MSLE, if
$b_n\asymp n^{-1/5}$:
%
\begin{eqnarray}\label{CLCSMSLE}
&& \sqrt{n} \biggl\{\tilde F_n^{\MSLE}(t)-F_0(t)-
\frac12b_n^2 \biggl\{ f_0'(t)+
\frac{2f_0(t)g'(t)}{g(t)} \biggr\}\int u^2K(u) \,du \biggr\} \Bigm/\sigma_n\hspace*{-25pt}
\nonumber\\[-8pt]\\[-8pt]\nonumber
&&\qquad \stackrel{{\mathcal D}}\longrightarrow N (0,1 ),
\end{eqnarray}
as $n\to\infty$, where $\sigma_n$ is defined as in (\ref
{CLCSSMLE}). Note that (\ref{CLCSSMLE}) and (\ref
{CLCSMSLE}) only differ in the bias term $b_n^2f_0(t)g'(t)/g(t)$.

\subsection{The SMLE and MSLE for the interval censoring, case 2, model}\label{subsectionintervalcens2model}
For interval censoring, case 2, we cannot rely on explicit
representations, as in the current status model. For the SMLE, we only
have a representation of type (\ref{SMLCS}) via the solution $\phi
^{\IC}$ of an integral equation, and we have to follow arguments
analogous to the arguments in \cite{piet96,GeGr96} and
\cite{GeGr97}.

I\vspace*{1pt}n the separated case (specified by Condition \ref{condition1} below),
the integral equation (in $\phi=\phi^{\IC}$) is given by
%
\begin{eqnarray}\label{F0-Fn-intequationsep}
\phi(u)&=&d_F(u) \biggl\{b^{-1}K
\bigl((t-u)/b \bigr)+\int_{v>u}\frac
{\phi(v)-\phi(u)}{F(v)-F(u)} g(u,v) \,dv
\nonumber\\[-8pt]\\[-8pt]\nonumber
&&\hspace*{109pt}{} -
\int_{v<u} \frac{\phi(u)-\phi(v)}{F(u)-F(v)} g(v,u) \,dv \biggr\},
\end{eqnarray}
where we take either $F=\hat F_n$ or $F=F_0$, and where
%
\begin{equation}
\label{dF} d_F(u)=\frac{F(u)\{1-F(u)\}}{g_1(u)\{1-F(u)\}+g_2(u)F(u)}.
\end{equation}

Moreover, let the function $\theta^{\IC}_{t,b,F}$ be defined by
%
\begin{eqnarray}\label{deftheta}
&& \theta^{\IC}_{t,b,F}(u,v,\delta_1,
\delta_2)
\nonumber\\[-8pt]\\[-8pt]\nonumber
&&\qquad =-\frac{\delta
_1\phi_{t,b,F}^{\IC}(u)}{F(u)}-\frac{\delta_2\{\phi
^{\IC}_{t,b,F}(v)-\phi^{\IC}_{t,b,F}(u)\}}{F(v)-F(u)}+ \frac{\delta
_3\phi^{\IC}_{t,b,F}(v)}{1-F(v)},\nonumber
\end{eqnarray}
where $u<v$. Then, as in \cite{GeGr97}, we have the representation
\begin{eqnarray*}
&& \int\IK\bigl((t-u)/b \bigr) \,d (\hat F_n-F_0 ) (u)
\\
&&\qquad =
\int\theta^{\IC}_{t,b,\hat F_n}(u,v,\delta_1,\delta
_2) \,dP_0(u,v,\delta_1 d_2)
\\
&&\qquad =\int\frac{\phi^{\IC}_{t,b,\hat F_n}(u)}{\hat F_n(u)}F_0(u)g_1(u) \,du
-\int
\frac{\phi^{\IC}_{t,b,\hat F_n}(v)}{1-\hat F_n(v)}\bigl\{1-F_0(v)\bigr
\} g_2(v) \,dv
\\
&&\quad\qquad{} +\int\frac{\phi^{\IC}_{t,b,\hat F_n}(v)-\phi
^{\IC}_{t,b,\hat F_n}(u)}{\hat F_n(v)-\hat F_n(u)}\bigl\{F_0(v)-F_0(u)
\bigr\} g(u,v) \,du \,dv.
\end{eqnarray*}
Using the theory in \cite{GeGr97} again, we get
that $\phi^{\IC}_{t,b,F_0}$ gives as an approximation for
$n \operatorname{var}(\tilde F_n(t))$:
\begin{eqnarray*}
&&E \theta^{\IC}_{t,b,F_0}(T_1,U_1,
\Delta_{11},\Delta_{12})^2
\\
&&\qquad =\int\frac{\phi^{\IC}_{t,b,F_0}(u)^2}{F_0(u)}g_1(u) \,du +\int\frac{ \{
\phi^{\IC}_{t,b,F_0}(v)-\phi
^{\IC}_{t,b,F_0}(u) \}^2}{F_0(v)-F_0(u)}h(u,v) \,du
\,dv
\\
&&\quad\qquad{}
+\int\frac
{\phi^{\IC}_{t,b,F_0}(v)^2}{1-F_0(v)}g_2(v) \,dv.
\end{eqnarray*}

Taking $b_n\asymp n^{-1/5}$ and defining
\[
\sigma_n^2=E \theta^{\IC}_{t,b_n,F_0}(T_1,U_1,
\Delta_{11},\Delta_{12})^2, %
\]
we get
%
\begin{eqnarray}\label{limitvarSMLE}
\qquad \lim_{b_n\downarrow0}b_n\sigma
_n^2
&=& d_{F_0}(t) \biggl\{1+d_{F_0}(t)\int_{v>t}
\frac{g(t,v)}{F_0(v)-F_0(t)} \,dv
\nonumber\\[-8pt]\\[-8pt]\nonumber
&&\hspace*{39pt}{} +d_{F_0}(t)\int_{v<t}
\frac{g(v,t)}{F_0(t)-F_0(v)} \,dv \biggr\} ^{-1}\int K(u)^2 \,du,
\end{eqnarray}
where $d_{F_0}$ is defined by (\ref{dF}).
This means, as we shall show below, that the limit variance for the
SMLE is again (as in the current status model) equal to the limit
variance of the MSLE. This leads to Conjecture 11.15 in \cite{pietgeurt14}:
%
\begin{eqnarray}\label{CLSSMLEIC}
\qquad\sqrt{n} \biggl\{\tilde F_n^{\SML}(t)-F_0(t)-
\frac12b_n^2f_0'(t)\int
u^2K(u) \,du \biggr\}\Bigm/\sigma_n \stackrel{{\mathcal D}}
\longrightarrow N(0,1),
\nonumber\\[-8pt]\\[-12pt]
\eqntext{n\to\infty,}
\end{eqnarray}
under the conditions given in \cite{pietgeurt14}. This also means
that the asymptotic bias is of the same form as for the SMLE in the
current status model, which is much simpler than the bias of the MSLE.

Throughout this paper, we will assume that the following conditions are
satisfied, which were also assumed in \cite{GeGr96} and \cite{GeGr97}.

%
\begin{condition}\label{condition1}
\textup{(S1)} $g_1$ and $g_2$ are continuous, with $g_1(x)+g_2(x) > 0$
for all $x \in[0,M]$.

\textup{(S2)} ${\mathbb P}\{V-U < \varepsilon\}=0$ for some
$\varepsilon$ with $0 < \varepsilon\leq1/2
M$, so $g$ does not have mass close to the diagonal; this is called
the separated case.

\textup{(S3)} $(u,v)\mapsto g(u,v)$ is continuous on $\{(x,y)\dvtx 0\le
x<y<M\}$ and is zero outside this set. Moreover, $g(u,v)=0$ if
$v-u<\varepsilon$.

\textup{(S4)} $F$ is a continuous distribution function with support
$[0,M]$; F satisfies
\[
F(u)-F(t) \geq c>0,\qquad\mbox{if }u-t\ge\varepsilon. %
\]

\textup{(S5)} The partial derivatives $\partial_1g(t,u)$ and
$\partial_2g(t,u)$ exist, except for at most a
countable number of points, where left and right derivatives exist. The
derivatives are bounded, uniformly over $t$ and $u$.

\textup{(S6)} If both $G_1$ and $G_2$ put zero mass on some set
$A$, then $F$ has zero mass on
$A$ as well, so $F \ll H_1+H_2$. This means that $F$ does not have mass
on sets in which no
observations can occur.
\end{condition}

Note that (S1) implies that $d_F$, defined by (\ref{dF}), is bounded.
Conditions (S2) and (S4) are needed to avoid singularity in the
integral equation: if $F(x)-F(t)$ becomes very small, we have
$g(t,x)=0$. A picture of an observation density, satisfying the above
conditions, is shown in Figure~\ref{figg}; $g$ is defined by
%
\begin{equation}
\label{defg} g(x,y)=6(y-x-\varepsilon)^2/\bigl\{(2-x-\varepsilon
) (2-\varepsilon)\bigr\} ^2,\qquad x+\varepsilon<y,
\end{equation}
on the triangle with vertices $(0,\varepsilon)$, $(0,2)$ and
$(2-\varepsilon,2)$, where $\varepsilon=0.1$.

We use the following conditions for the kernel estimators.
%
\begin{condition}[(Conditions on the kernel estimators)]\label{condition2}
We assume that $\tilde h_{nj}$ and $\tilde h_n$ are kernel estimators
of $h_{0j}$ and $h_0$, respectively, defined by (\ref{defhnj}) and
(\ref{defhn}), respectively, for a symmetric continuously
differentiable kernel $K$ of type (\ref{defK}), with compact support.
Moreover, for points near the boundary, boundary kernels are used, with
coefficients $\alpha(t)$ and $\beta(t)$, defined by (\ref{coeff1})
and (\ref{coeff2}), respectively, where the functions $\alpha$,
$\beta$, and its derivatives $\alpha'$ and $\beta'$ are bounded
on $[0,1]$. We assume:
%
\begin{equation}
\label{deflowerend} 0=\inf\biggl\{t\in[0,M]\dvtx \tilde h_{n1}(t)\vee\int
_{u=0}^t\tilde h_n(u,t) \,du>0 \biggr\}
\end{equation}
and
%
\begin{equation}
\label{defupperend} M=\sup\biggl\{t\in[0,M]\dvtx \tilde h_{n2}(t)\vee\int
_{u=t}^M\tilde h_n(t,u) \,du>0 \biggr\}.
\end{equation}
\end{condition}

An example of a kernel estimate, satisfying the conditions of Condition
\ref{condition2}, is given by kernel estimates which use the triweight
kernel, defined by (\ref{defK}). For this kernel, the weight
functions $\alpha$ and $\beta$, used in constructing the boundary
kernel, are decreasing on $[0,1]$, and the derivatives are bounded on $[0,1]$.
Using Condition \ref{condition1}, we give a characterization in terms
of necessary and sufficient (duality) conditions for the MSLE in
Section~\ref{sectionconsistency}. In that section, we also prove
consistency of the MSLE, using techniques, similar to the method, used
in \cite{GrWe92}, Part II, Section~4.

In Section~\ref{sectionalgorithms}, we discuss algorithms for
computing the MSLE: the EM algorithm and an iterative convex minorant
algorithm. The iterative convex minorant algorithm is an adapted
version of the algorithm, introduced in \cite{Gr91} and (again in)~\cite{GrWe92}. It turns out that the latter algorithm performs best
in our experiments. The EM algorithm is very slow and, therefore, not
suitable for larger sample sizes or simulation purposes.

In Section~\ref{sectionnormality}, we will prove asymptotic normality
of the MSLE at a fixed interior point of the domain of definition
(Theorem \ref{thasympnorsep}).
In this paper, we concentrate on the ``separated case,'' where
$U_i-T_i\ge\varepsilon$ for some $\varepsilon>0$, as in \cite
{GeGr96} and \cite{GeGr97}. This case seems to be the most important
case, and also to be the usual situation in medical statistics. The
nonseparated case is rather different and has its own specific
difficulties. The behavior of the MLE and SMLE in this situation is
discussed in \cite{GeGr99,piet96} and \cite{piettom11},
but the theory is still rather incomplete, even for the MLE. There is a
conjecture for its asymptotic distribution, put forward in \cite
{Gr91} and \cite{GrWe92}, but this conjecture has not been proved up
till now, although a simulation study, supporting the conjecture is
given in \cite{piettom11}.
The theory for the MSLE in this situation has still not been developed.

\section{Characterization of the MSLE and consistency}\label{sectionconsistency}

Let, for an estimate $\tilde h_n$ of $h_0$, satisfying
%
\begin{equation}
\label{separationconditionestimate} \tilde h_n(t,u)=0,\qquad u-t<\varepsilon,
\end{equation}
for some $\varepsilon>0$, the nabla function $\nabla_F$ be defined by
%
\begin{eqnarray}\label{nablaIC}
\nabla_F(u) &=& \frac{\tilde h_{n1}(u)}{F(u)}-\frac{\tilde
h_{n2}(u)}{1-F(u)} +
\int_{v=0}^u \frac{\tilde h_{n}(v,u)}{F(u)-F(v)} \,dv
\nonumber\\[-8pt]\\[-8pt]\nonumber
&&{} -\int
_{v=u}^M \frac{\tilde h_{n}(u,v)}{F(v)-F(u)} \,dv,\nonumber
\end{eqnarray}
if $0<F(u)<1$. If $F(u)=0$ or $F(u)=1$, we define $\nabla_F(u)=0$.

Then, similarly to the ordinary MLE, the MSLE can be characterized by
the so-called Fenchel duality conditions.

%
\begin{lemma}
\label{lemmacharacterization}
Let $\tilde h_n$ satisfy (\ref{separationconditionestimate}), for
some $\varepsilon>0$. Then the distribution\vspace*{2pt} function $\hat F_n$
maximizes (\ref{criterionfunction1}) if and only if $\hat F_n$ is
continuous on $[0,M]$ and satisfies
the conditions
%
\begin{equation}
\label{fenchel1} \int_{v=t}^M
\nabla_{\hat F_n}(v) \,dv\le0,\qquad t\in[0,M)
\end{equation}
and
%
\begin{equation}
\label{fenchel2} \int_0^M\nabla_{\hat F_n}(v)
\hat F_n(v) \,dv=0.
\end{equation}
Moreover, if $t\in[0,M)$ is a point of increase of $\hat F_n$, that is,
%
\begin{equation}\label{increasepoint}
\cases{ \hat F_n(u)-\hat F_n \bigl(u'\bigr)>0,
\cr
\qquad \mbox{for all $u,u'\in[0,M]$ such that $u'<t<u$,} &\quad if $t>0$,
\cr
\hat F_n(u)>0,
\cr
\qquad \mbox{$u\in(0,M]$,} &\quad if $t=0$,}
\end{equation}
we have
%
\begin{equation}
\label{fenchel3} \nabla_{\hat F_n}(t)=0\quad\mbox{and}\quad\int
_t^M\nabla_{\hat F_n}(v) \,dv=0.
\end{equation}
\end{lemma}

The proof of this lemma is given in the \hyperref[app]{Appendix}.

Note that if $\nabla_{F}(t)=0$ for all $t\in(0,M)$, where $\nabla
_{F}$ is defined by (\ref{nablaIC}), the conditions of Lemma \ref
{lemmacharacterization} are satisfied for $F$, and hence $F$ would be
the MSLE if it also would be a distribution function. But
unfortunately, the function $F$ satisfying $\nabla_{F}(t)=0$ for all
$t\in(0,M)$ need not be monotone. We will call a function~$\tilde
F_n$, satisfying $\nabla_{\tilde F_n}(t)=0$, $t\in(0,M)$, a \textit{plug-in estimator} or
\textit{naive estimator} (as\vspace*{1.5pt} in~\cite
{pietgeurtbirgit10}). This plug-in estimator will be further studied
in Section~\ref{sectionnormality} in the proof of the local
asymptotic normality of the MSLE, where it will be shown that the MSLE
is indeed locally asymptotically equivalent to this plug-in estimator.

%
\begin{corollary}
\label{corfenchel}
Let $\tilde h_n$ satisfy (\ref{separationconditionestimate}), for
some $\varepsilon>0$. Then the distribution function $\hat F_n$
maximizes (\ref{criterionfunction1}) if and only if $\hat F_n$ is
continuous on $[0,M)$, $\hat F_n(M)>0$, and if $\hat F_n$ satisfies the
conditions
%
\begin{equation}
\label{fenchel1a} \int_0^t\nabla_{\hat F_n}(v)
\,dv\ge0,\qquad t\in(0,M)
\end{equation}
and
%
\begin{equation}
\label{fenchel2a} \int_0^M\nabla_{\hat F_n}(v)
\,dv=0.
\end{equation}
Moreover, if $t\in[0,M)$ is a point of increase of $\hat F_n$, that is,
satisfies condition (\ref{increasepoint}) of Lemma \ref
{lemmacharacterization}, then
%
\begin{equation}
\label{fenchel3a} \nabla_{\hat F_n}(t)=0\quad\mbox{and}\quad\int
_0^t\nabla_{\hat F_n}(v) \,dv=0.
\end{equation}
\end{corollary}

\begin{pf}
Suppose $\hat F_n$ maximizes $\ell(F)$. Defining
\[
F_{\delta}(t)=\bigl\{1-(1+\delta) \bigl(1-\hat F_n(t)\bigr)
\bigr\}\vee0,\qquad t\in[0,M], %
\]
we find:
%
\begin{equation}
\label{1-F-condition} \lim_{\delta\to0}\frac{\ell(F_{\delta})-\ell
(\hat
F_n)}{\delta}= -\int
_0^M \nabla_{\hat F_n}(u) \bigl\{1-\hat
F_n(u) \bigr\} \,du=0.
\end{equation}
So if $\hat F_n$ maximizes $\ell(F)$, (\ref{fenchel2a}) follows from
(\ref{1-F-condition})
and (\ref{fenchel2}) of Lemma \ref{lemmacharacterization}.
\[
\int_0^M \nabla_{\hat F_n}(u) \,du=0.
\]
This implies
\[
\int_0^t\nabla_{\hat F_n}(v) \,dv=-\int
_t^M\nabla_{\hat F_n}(v) \,dv, %
\]
and condition (\ref{fenchel3a}) now also follows.

Conversely, if the conditions of the corollary hold, we get
\begin{eqnarray*}
&& \int_0^M\hat F_n(u)
\nabla_{\hat F_n}(u) \,du
\\
&&\qquad =\hat F_n(M)\int_0^M
\nabla_{\hat F_n}(v) \,dv+\int_{t=0}^M\int
_{v=0}^u\nabla_{\hat
F_n}(v) \,dv \,d\hat
F_n(u)
\\
&&\qquad =\hat F_n(M)\int_0^M
\nabla_{\hat F_n}(v) \,dv=0,
\end{eqnarray*}
implying condition (\ref{fenchel2}) of Lemma \ref
{lemmacharacterization}. The other conditions of Lemma \ref
{lemmacharacterization} follow similarly.
\end{pf}

We now simplify the conditions somewhat, in view of the iterative
convex minorant algorithm, to be discussed in Section~\ref{sectionalgorithms}. Multiplying $\nabla_F$ by $F(1-F)$ yields the function
%
\begin{eqnarray}
\label{nablaIC2} \overline{\nabla}_F(u)&=&\tilde h_{n1}(u)
\bigl\{1-F(u)\bigr\}-\tilde h_{n2}(u)F(u)
\nonumber
\\
&&{}+F(u)\bigl\{1-F(u)\bigr\}
\\
&&\hspace*{10pt}{}\times \biggl\{\int_{v=0}^u
\frac{\tilde
h_{n}(v,u)}{F(u)-F(v)} \,dv
-\int_{v=u}^M
\frac{\tilde h_{n}(u,v)}{F(v)-F(u)} \,dv \biggr\}.\nonumber
\end{eqnarray}

%
\begin{corollary}
\label{corfenchel2}
Let $\tilde h_n$ satisfy (\ref{separationconditionestimate}), for
some $\varepsilon>0$ and let the function $\overline{\nabla}_F$ be
defined by (\ref{nablaIC2}). Then the distribution function $\hat
F_n$ maximizes (\ref{criterionfunction1}) if and only if $\hat
F_n(M)>0$, and $\hat F_n$ is continuous on $[0,M)$ and satisfies the conditions
%
\begin{equation}
\label{fenchel1b} \int_0^t\overline{
\nabla}_{\hat F_n}(v) \,dv\ge0,\qquad t\in[0,M]
\end{equation}
and
%
\begin{equation}
\label{fenchel2b} \int_0^M\nabla_{\hat F_n}(v)
\,dv=0.
\end{equation}
Moreover, if $t\in[0,M)$ is a point of increase of $\hat F_n$, that
is, satisfies condition (\ref{increasepoint}) of Lemma \ref
{lemmacharacterization}, then
%
\begin{equation}
\label{fenchel3b} \overline{\nabla}_{\hat F_n}(t)=0\quad\mbox{and}\quad
\int_0^t\overline{\nabla}_{\hat F_n}(v)
\,dv=0.
\end{equation}
\end{corollary}

\begin{pf}
We have, for $t\in(0,M)$,
\begin{eqnarray*}
\int_0^t\overline{\nabla}_{\hat F_n}(v)
\,dv&=&\int_0^t \hat F_n(v)\bigl\{ 1-
\hat F_n(v)\bigr\}\nabla_{\hat F_n}(v) \,dv.
\end{eqnarray*}
Furthermore,
\begin{eqnarray*}
&&\int_a^t \hat F_n(v)\bigl\{1-
\hat F_n(v)\bigr\}\nabla_{\hat F_n}(v) \,dv
\\
&&\qquad =\hat F_n(t)\bigl\{1-\hat F_n(t)\bigr\}\int
_{u=0}^t \nabla_{\hat F_n}(u) \,du
\\
&&\quad\qquad{} -\int
_0^t \bigl\{1-2\hat F_n(u)\bigr\}
\int_{v=0}^u \nabla_{\hat F_n}(v) \,dv \,d\hat
F_n(u)
\\
&&\qquad =\hat F_n(t)\bigl\{1-\hat F_n(t)\bigr\}\int
_{u=0}^t \nabla_{\hat F_n}(u) \,du.
\end{eqnarray*}
Hence, condition (\ref{fenchel1b}) is equivalent to condition (\ref
{fenchel1a}) of Corollary \ref{corfenchel}.
Relation~(\ref{fenchel3b}) follows similarly, and (\ref{fenchel2b})
is the same as (\ref{fenchel2a}).
\end{pf}

The preceding results imply the consistency of the MSLE. The proof,
which is given in the \hyperref[app]{Appendix}, is somewhat analogous to the proof of
the consistency of the MLE in \cite{GrWe92}.

%
\begin{theorem}[(Consistency of the MSLE)]\label{thconsistency}
Let Condition \ref{condition1} be satisfied on~$[0,M]$ for the\vspace*{1pt}
distribution function $F_0$ and the observation density $g$. Moreover,
let $\tilde h_{nj}$ and $\tilde h_n$ be kernel estimators of $h_{0j}$
and $h_0$, respectively, of the type defined in Condition \ref
{condition2}. Finally, let $\hat F_n$ be the MSLE of $F_0$. Then, with
probability one,
\[
\lim_{n\to\infty}\hat F_n(t)=F_0(t),
\]
for each $t\in[0,M)$. The convergence is uniform on each subinterval
$[a,b]$ of~$(0,M)$.
\end{theorem}

The proof of this theorem is given in the \hyperref[app]{Appendix}.

\section{Algorithms}\label{sectionalgorithms}
We explained in Section~\ref{sectionintro} that the MSLE can be
computed for current status data via a continuous cusum diagram. In the
present case we do not have a similar algorithm, which computes the
MSLE in one step.
The EM algorithm is based on the following ``self-consistency equations''
\begin{eqnarray*}
\hat f_n(t) &=& \biggl\{\int_t^M
\frac{\tilde h_{n1}(v)}{\hat F_n(v)} \,dv +\int_0^t
\frac{\tilde h_{n2}(v)}{1-\hat F_n(v)} \,dv
\\
&&\hspace*{19pt}{} +\int_{v<t<u} \frac{\tilde
h_n(v,u)}{\hat F_n(u)-\hat F_n(v)} \,dv \,du
\biggr\}\hat f_n(t),
\end{eqnarray*}
where $\hat f_n(t)=\hat F_n'(t)$. This yields the iteration steps
%
\begin{eqnarray}
\label{EM-iteration} f^{(k+1)}(t)&=& \biggl\{\int_t^M
\frac{\tilde h_{n1}(v)}{F^{(k)}(v)} \,dv +\int_0^t
\frac{\tilde h_{n2}(v)}{1-F^{(k)}(v)} \,dv
\nonumber\\[-8pt]\\[-8pt]\nonumber
&&\hspace*{16pt}{} +\int_{v<t<u} \frac{\tilde
h_n(v,u)}{F^{(k)}(u)-F^{(k)}(v)} \,dv \,du
\biggr\}f^{(k)}(t).
\end{eqnarray}
One can indeed use a discretized version of (\ref{EM-iteration}) to
compute the MSLE, but the EM algorithm is (as is usual for this type of
problem with many parameters) very slow. Simply enhancing the EM
algorithm by a Newton step is also not helpful because of the many
constraints the solution has to satisfy, leading to very small
``feasible steps.'' For this reason, a Newton-improved EM algorithm
does not improve very much on the EM algorithm itself.

In our experience, the fastest algorithm is a combination of the EM
algorithm with a version of the iterative convex minorant (ICM)
algorithm, introduced in \cite{Gr91} and \cite{GrWe92}. We use a
sequence of cusum diagrams
%
\begin{equation}
\label{cusumsmoothcurstat2} \bigl(W_n^{(k)}(t),V_n^{(k)}(t)
\bigr),\qquad t\in[0,M], k=0,1,2,\ldots,
\end{equation}
for which we compute the greatest convex minorants at each $k$th step.
We alternate this with an EM-step (the combination is sometimes called
the ``hybrid algorithm''). The cumulative weight function $W_n^{(k)}$
is of the form
\[
W_n^{(k)}=\int_0^t
w_n^{(k)}(u) \,du,\qquad t\ge0, %
\]
for suitably (but somewhat arbitrarily) chosen weights $w_n^{(k)}$, and
the cusum function $V_n^{(k)}$ is of the form:
\[
V_n^{(k)}(t)=\int_0^t
F^{(k)}(u)w_n^{(k)}(u) \,du+\int_0^t
\overline{\nabla}_{F^{(k)}}(u) \,du,\qquad t\ge0, %
\]
where, for a distribution function $F$, $\overline{\nabla}_F$ is the
function, defined by (\ref{nablaIC2}),
evaluated at $F=F^{(k)}$. The idea is that the iterations force the
Fenchel duality conditions (\ref{fenchel1b}) and (\ref{fenchel2b})
to be satisfied at the end of the iterations.

The following weight function, chosen by taking the diagonal elements
of the Hessian matrix, corresponding to the function $\overline{\nabla
}_{F}$, gave good convergence results in our simulation study of the MSLE:
\begin{eqnarray*}
w_n^{(k)}(t)&=&\tilde h_{n1}(t)+\tilde
h_{n2}(t)
\\
&&{}- \bigl\{1-2F_n^{(k)}(t) \bigr\}
\\
&&\hspace*{10pt}{}\times \biggl\{\int
_{u=0}^t\frac
{\tilde h_n(u,t)}{F^{(k)}(t)-F^{(k)}(u)} \,du -\int
_{u=t}^M\frac{\tilde h_n(t,u)}{F^{(k)}(u)-F^{(k)}(t)} \,du \biggr\}
\\
&&{}+F^{(k)}(t)\bigl\{1-F^{(k)}(t)\bigr\}\int
_{u=0}^t\frac{\tilde
h_n(u,t)}{\{F^{(k)}(t)-F^{(k)}(u)\}^2} \,du
\\
&&{} +\int
_{u=t}^M\frac{\tilde h_n(t,u)}{\{F^{(k)}(u)-F^{(k)}(t)\}^2} \,du.
\end{eqnarray*}
To prevent divergence of the algorithm, Armijo's line search method, as
implemented in \cite{geurt98}, was used for determining the step size
at each iteration. The integrals were computed by a discrete
approximation, using Riemann sums.

Note that, in the case of current status data, the function $\overline
{\nabla}_F$ is just given by
%
\begin{equation}
\label{nablaCS} \overline{\nabla}_{F}=\tilde h_n(t)-
\tilde g_n(t)F(t),
\end{equation}
from which we can compute $\hat F_n$ in one step.

\section{Asymptotic distribution}\label{sectionnormality}

\subsection{Main result and road map}
\label{subsectionroadmap}

We will prove the following theorem.

%
\begin{theorem}
\label{thasympnorsep}
Let condition (\ref{condition1}) be satisfied. Moreover,
let $F_0$ be twice differentiable, with a bounded continuous
derivative\vadjust{\goodbreak}
$f_0$ on the interior of $[0,M]$, which is bounded away from zero on
$[0,M]$, with a finite positive right limit at $0$ and a positive left
limit at $M$. Also, let $f_0$ have a bounded continuous\vadjust{\goodbreak} derivative on
$(0,M)$ and let $g_1$ and $g_2$ be twice differentiable on the interior
of their supports $S_1$ and $S_2$, respectively. Furthermore, let the
joint density $g$ of the pair of observation times $(T_i,U_i)$ have a
bounded (total) second derivative on $\{(x,y)\dvtx 0<x<y<M\}$. Suppose that
$X_i$ is independent of $(T_i,U_i)$, and let $d_{F_0}$ be defined by
\[
d_{F_0}(v)=\frac{F_0(v)\{1-F_0(v)\}}{g_1(v)\{1-F_0(v)\}+F_0(v)g_2(v)}. %
\]
Then, if $b_n\asymp n^{-1/5}$, we have for each $v\in(0,M)$,
\[
\sqrt{nb_n} \bigl\{\hat F_n(v)-F_0(v)-\beta
(v)b_n^2 \bigr\} \stackrel{{\mathcal D}}\longrightarrow N
\bigl(0,\sigma(v)^2 \bigr), %
\]
where $N (0,\sigma(v)^2 )$ is a normal distribution with
first moment zero and variance $\sigma(v)^2$, and where, defining
%
\begin{equation}\label{sigma1}
\sigma_1(v)=1+d_{F_0}(v) \biggl\{\int
_{t<v}\frac
{g(t,v)}{F_0(v)-F_0(t)} \,dt +\int_{w>v}
\frac{g(v,w)}{F_0(w)-F_0(v)} \,dw \biggr\},\hspace*{-25pt}
\end{equation}
the variance $\sigma(v)^2$ is given by
%
\begin{equation}
\label{as-variance} \sigma(v)^2=\frac{d_{F_0}(v)}{
\sigma_1(v)}\int
K(u)^2 \,du.
\end{equation}
Defining
%
\begin{eqnarray}
\label{bias1} \beta_1(v)&=&\frac{1}{2\sigma_1(v)} \biggl\{
\frac{\{1-F_0(v)\}
h_1''(v)-F_0(v)h_2''(v)}{g_1(v)\{1-F_0(v)\}+F_0(v)g_2(v)}
\nonumber
\\
&&\hspace*{37pt}{} +d_{F_0}(v) \biggl\{\int_{t=0}^v
\frac
{(\fraca{\partial^2}{\partial v^2})h_0(t,v)}{F_0(v)-F_0(t)} \,dt
\\
&&\hspace*{82pt}{}-\int_{u=v}^M
\frac{(\fraca{\partial^2}{\partial v^2})h_0(v,u)}{F_0(u)-F_0(v)} \,du \biggr\} \biggr\} \int u^2K(u) \,du,\nonumber
\end{eqnarray}
the bias $\beta(v)$ is given by
\begin{eqnarray*}
\beta(v)&=&\beta_1(v)+\frac{d_{F_0}(v)}{\sigma_1(v)} \biggl\{\int
_{u=0}^v\frac{g(u,v)\beta_1(u)}{F_0(v)-F_0(u)} \,du +\int
_{u=v}^M\frac{g(v,u)\beta_1(u)}{F_0(u)-F_0(v)} \,du \biggr\}.
\end{eqnarray*}
\end{theorem}

%
\begin{remark}
The asymptotic bias of the MSLE is of a very complicated form,
certainly compared to the asymptotic bias of the SMLE, which is just
\[
\frac12f_0'(t)b_n^2\int
u^2 K(u) \,du; %
\]
see (\ref{CLSSMLEIC}). It would be nice if some simplification could
be found. Note however, that also in the current status model the
asymptotic bias of the MSLE is more complicated than the asymptotic
bias of the SMLE, since the derivatives of the estimates of the
observation density come into play.
\end{remark}

We now first give a ``road map'' of the proof of Theorem \ref{thasympnorsep}.
Our starting point is given by the duality conditions (\ref
{fenchel1b}) and (\ref{fenchel2b}).
It is clear that if we would\vadjust{\goodbreak} have equality in (\ref{fenchel1b})
instead of inequality, we would get the following relation by
differentiating w.r.t. $t$:
%
\begin{eqnarray}
\label{nabla=0} \overline{\nabla}_F(t)&=&\tilde h_{n1}(t)
\bigl\{1-F(t)\bigr\}-\tilde h_{n2}(t)F(t)
\nonumber
\\
&&{} +F(t)\bigl\{1-F(t)\bigr\} \biggl\{\int_{v=0}^t
\frac{\tilde
h_{n}(v,t)}{F(t)-F(v)} \,dv -\int_{v=t}^M
\frac{\tilde h_{n}(t,v)}{F(v)-F(t)} \,dv \biggr\}\hspace*{-15pt}
\\
&=& 0.\nonumber
\end{eqnarray}
Conversely, if $F$ solves (\ref{nabla=0}) for each $t\in(0,M)$ and
$F$ is a distribution function such that $F(t)\in(0,1)$, for each
$t\in(0,M)$, $F$ also satisfies (\ref{fenchel1b}) and (\ref
{fenchel2b}) and is therefore the MSLE.

The solution of equation (\ref{nabla=0}) takes the role of the plugin
estimator (\ref{pluginCS}) in the current status model. In the proof
of the central limit theorem for the MSLE for the current status model,
it was shown that the solution (in $F$) of (\ref{nabla=0}) is a
distribution function on a subinterval $(a,b)$ of $[0,M]$ for large $n$
with high probability, where we can take $a>0$ arbitrarily close to $0$
and $b<M$ arbitrarily close to $M$. In the present case, we prove the
stronger fact that the solution\vspace*{1pt} of (\ref{nabla=0}) is a
(sub)distribution function on $[0,M]$ itself. This implies that the
MSLE $\hat F_n$ coincides with the solution of (\ref{nabla=0}) on the
interval $(0,M)$ for large $n$ with high probability.

To show that the solution (in $F$) of (\ref{nabla=0}) is with high
probability a (sub)distribu\-tion function on $[0,M]$ for large $n$, we
first show that the solution is close to the solution of the \textit{linear} integral equation
%
\begin{eqnarray}
\label{linearasympequation1} && F(t)-F_0(t) +d_{F_0}(t) \biggl\{\int
_{u=0}^t\frac{g(u,t)\{F(t)-F_0(t)-F(u)+F_0(u)\}
}{F_0(t)-F_0(u)} \,du
\nonumber
\\
&&\hspace*{90pt}\quad{}-\int
_{u=t}^M\frac{g(t,u)\{F(u)-F_0(u)-F(t)+F_0(t)\}}{F_0(u)-F_0(t)} \,du
\biggr\}
\nonumber\\[-8pt]\\[-8pt]\nonumber
&&\qquad =\frac{\tilde h_{n1}(t)\{1-F_0(t)\}-\tilde h_{n2}(t)F_0(t)}{\{
1-F_0(t)\}g_1(t)+F_0(t)g_2(t)}
\nonumber
\\
&&\quad\qquad{} +d_{F_0}(t) \biggl\{\int_{u<t}
\frac
{\tilde h_n(u,t)}{F_0(t)-F_0(u)} \,du-\int_{u>v}\frac{\tilde
h_n(t,u)}{F_0(u)-F_0(t)} \,du \biggr
\},\nonumber
\end{eqnarray}
where $d_{F_0}$ is defined by (\ref{dF}), with $F=F_0$.

We next show that the ``toy estimator,'' solving the equation
%
\begin{eqnarray}
\label{linearasympequation2} && \bigl\{F(t)-F_0(t) \bigr\} \biggl\{1+
d_{F_0}(t) \biggl\{\int_{u=0}^t
\frac{g(u,t)}{F_0(t)-F_0(u)} \,du\nonumber\\
&&\hspace*{130pt}+\int_{u=t}^M
\frac{g(t,u)}{F_0(u)-F_0(t)} \,du \biggr\} \biggr\}\nonumber\hspace*{-5pt}
\\[-8pt]\\[-8pt]
&&\qquad =\frac{\tilde h_{n1}(t)\{1-F_0(t)\}-\tilde h_{n2}(t)F_0(t)}{\{
1-F_0(t)\}g_1(t)+F_0(t)g_2(t)}\nonumber
\\
&&\quad\qquad{} +d_{F_0}(t) \biggl\{\int_{u<t}
\frac{\tilde h_n(u,t)}{F_0(t)-F_0(u)} \,du-\int_{u>v}\frac{\tilde
h_n(t,u)}{F_0(u)-F_0(t)} \,du \biggr
\},\nonumber
\end{eqnarray}
where the off-diagonal terms
%
\begin{eqnarray}\label{off-diag}
&& -d_{F_0}(t) \biggl\{\int_{u=0}^t
\frac{g(u,t)\{F(u)-F_0(u)\}
}{F_0(t)-F_0(u)} \,du
\nonumber\\[-8pt]\\[-8pt]\nonumber
&&\hspace*{40pt}{}
+\int_{u=t}^M
\frac{g(t,u)\{F(u)-F_0(u)\}
}{F_0(u)-F_0(t)} \,du \biggr\}
\end{eqnarray}
on the left-hand side of (\ref{linearasympequation1}) are omitted,
also solves (\ref{linearasympequation1}) to the right order, apart
from a deterministic shift term.
This last step is somewhat similar to a part of the proof of the
asymptotic distribution of the MLE for interval censoring under the
separation condition in \cite{piet96}. However, in the latter case a
corresponding ``off-diagonal'' term (\ref{off-diag}) plays no role
asymptotically, since for the MLE the contribution to the bias is of
lower order.
In this way, we have reduced the proof to the asymptotic equivalence of
the MSLE with the solution of (\ref{linearasympequation2}) on $[0,M]$.

\begin{figure}

\includegraphics{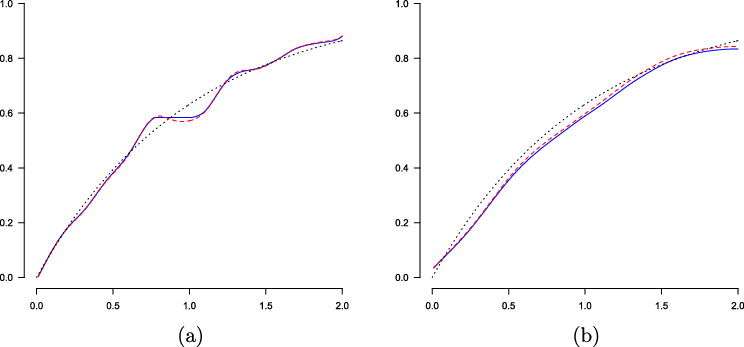}

\caption{\textup{(a)}~The MSLE (solid) and the toy estimator, solving equation
(\protect\ref{linearasympequation2}) (dashed), for a sample of size
$n=1000$ from the exponential distribution function $F_0(x)=1-\exp\{
-x\}$ (dotted), using bandwidth $b_n=n^{-1/5}\approx0.25119$. The
bivariate observation density $g$ is defined by (\protect\ref{defg}), where $\varepsilon=0.1$.
\textup{(b)}~The same, but now using the bandwidth $b_n=2n^{-1/5}\approx0.36411$.}\label{figMSLEinteq}
\end{figure}

A comparison of the MSLE and the toy estimator, solving (\ref
{linearasympequation2}), is shown in part (a) of Figure~\ref
{figMSLEinteq} for bandwidth $b_n=n^{-1/5}$. One can see that for
this bandwidth isotonization is still needed (the MSLE has derivative
zero on a piece in the middle of the interval). Note that the toy
estimator is not monotone, but has a very small distance to the MSLE.
If we take $b_n=2n^{-1/5}$, as in part (b) of Figure~\ref
{figMSLEinteq}, which seems a better choice in this case,
isotonization is not needed, except at the very last end of the
interval (where the MSLE has derivative zero). Note that, as $n\to
\infty$, the bandwidth will become smaller than $\varepsilon/2$,
where $\varepsilon$ is the separation distance in~(\ref{defg}), but
that this is still not the case in Figures~\ref{figMSLEinteq}.\vadjust{\goodbreak}

%

Replacing $\tilde h_{nj}$ by $h_{0j}$ in (\ref{nabla=0}), $j=1,2$, and
$\tilde h_n$ by $h_0$, we obtain the equation
\begin{eqnarray*}
&& h_{01}(t)\bigl\{1-F(t)\bigr\}-h_{02}(t)F(t)
\\
&&\qquad{} +F(t)\bigl
\{1-F(t)\bigr\} \biggl\{\int_{v=0}^t
\frac{h_0(v,t)}{F(t)-F(v)} \,dv -\int_{u=t}^M
\frac{h_0(t,u)}{F(u)-F(t)} \,du \biggr\}=0,
\end{eqnarray*}
which, using the definition of $h_{0j}$ and $h_0$, turns into
\begin{eqnarray*}
&& g_1(t)F_0(t)\bigl\{1-F(t)\bigr\}-g_2(t)
\bigl\{1-F_0(t)\bigr\}F(t)
\\
&&\quad{}
+F(t)\bigl\{1-F(t)\bigr\}
\\
&&\qquad{}\times \biggl\{\int_{v=0}^t
\frac{g(v,t)\{F_0(t)-F_0(u)\}
}{F(t)-F(v)} \,dv -\int_{u=t}^M
\frac{g(t,u)\{F_0(u)-F_0(t)\}}{F(u)-F(t)} \,du \biggr\}
\\
&&\qquad =0.
\end{eqnarray*}
This equation is clearly solved by $F_0$ itself.

This motivates us to consider the equation
\begin{eqnarray*}
\phi(t;h_1,h_2,h,F)&=&0,\qquad t\in[0,M],
\end{eqnarray*}
where
%
\begin{eqnarray}\label{fundamentaleq}
&& \phi(t;h_1,h_2,h,F)\nonumber
\\
&&\qquad = h_1(t)\bigl\{1-F(t)\bigr\}-h_2(t)F(t)
\\
&&\quad\qquad{}+F(t)\bigl\{1-F(t)\bigr\} \biggl\{\int_{v=0}^t
\frac{h(v,t)}{F(t)-F(v)} \,dv -\int_{u=t}^M
\frac{h(t,u)}{F(u)-F(t)} \,du \biggr\}.\nonumber\hspace*{-20pt}
\end{eqnarray}
The functions $h$ belong to a closed subset of the Banach space $C(S)$,
where $S$ is given by
\[
S=\bigl\{(x,y)\dvtx 0\le x\le x+\varepsilon_0\le y\le M\bigr\},
\]
for a fixed $\varepsilon_0>0$. We further define
\[
S_1=[0,M-\varepsilon_0],\qquad S_2=[\varepsilon
_0,M]. %
\]

We now define the space $E$ by
%
\begin{equation}
\label{def-E} E=C[S_1]\times C[S_2]\times C(S)\times
C[0,M],
\end{equation}
and put the following norm on $E$:
%
\begin{equation}
\label{fundamentalnorm} \bigl\llVert(h_1,h_2,h,F)\bigr\rrVert
=\max\bigl\{\llVert h_1\rrVert,\llVert h_2\rrVert,
\llVert h\rrVert,\llVert F\rrVert\bigr\},
\end{equation}
where the norms on the right-hand side denote the supremum norm, which
we also denote by $\llVert \cdot\rrVert $. Note that $E$ is a Banach
space for the
norm (\ref{fundamentalnorm}).

We will also need another norm on $C[S]$, defined by
%
\begin{equation}
\label{h-norm} \llVert h\rrVert_S=\sup_{t\in[0,M]}
\biggl\{\int_{u\dvtx (u,t)\in S}\bigl\llvert h(u,t)\bigr\rrvert \,du +\int
_{u\dvtx (t,u)\in S}\bigl\llvert h(t,u)\bigr\rrvert \,du \biggr\}.
\end{equation}
Note that this is indeed a norm on $C[S]$, since $\llVert h\rrVert
_S=0$ implies
$h=0$ and since the triangle inequality and homogeneity property for
scalars are obviously satisfied.\vadjust{\goodbreak}

The proof of Theorem \ref{thasympnorsep} now proceeds via a
sequence of lemmas. The proofs of these lemmas are given in the
\hyperref[app]{Appendix}. The first lemma tells us that we can apply the implicit
function theorem in Banach spaces to ensure that, locally, using the
norms just introduced, there is a unique solution $F$ to the equation
$\phi(t;h_1,h_2,h,F)=0$.

%
\begin{lemma}
\label{lemmaexistencesolution}
Let $F_0$, $h_{01}$, $h_{02}$ and $h_0$ satisfy the conditions of
Theorem \ref{thasympnorsep}. Furthermore, let the function $\phi$
be defined by (\ref{fundamentaleq}). Then there exists for all small
$\eta>0$ an open set $U$ in the Banach space $C[S_1]\times
C[S_2]\times C(S)$, endowed with the norm
\[
\bigl\llVert(h_1,h_2,h)\bigr\rrVert=\max\bigl\{
\llVert h_1\rrVert,\llVert h_2\rrVert,\llVert h\rrVert
\bigr\}, %
\]
such that, if $(h_1,h_2,h)\in U$, the equation
\[
\phi(t;h_1,h_2,h,F)=0,\qquad t\in[0,M],
\]
where $\phi$ is defined by (\ref{fundamentaleq}), has a unique
solution $F$ in the open ball $B(F_0,\eta)\subset C[0,M]$ with
midpoint $F_0$.
\end{lemma}

Having established the existence of a solution, we also consider the
derivative of the solution.

%
\begin{lemma}
\label{lemmaderivativesolution}
Let, under the conditions of Lemma \ref{lemmaexistencesolution}, for
a small $\eta>0$, $F\in B(F_0,\eta)$ be the solution of
\[
\phi(t;h_1,h_2,h,F)=0,\qquad t\in[0,M], %
\]
where $\phi$ is defined by (\ref{fundamentaleq}), and where $h_j$
has a bounded continuous derivative on the interior of $S_j$, having
finite limits approaching the boundary of $S_j$, for $j=1,2$.
Similarly, we suppose that $h$ is differentiable on the interior of its
support $S$ and has finite limits approaching the boundary of $S$.
Then, if $(h_1,h_2,h)\in U_{\delta}$, where $U_{\delta}$ is defined
by (\ref{U-delta}), the solution $F$ has a continuous and bounded
derivative for sufficiently small $\eta$ and $\delta$.
\end{lemma}

The following lemma will be used to show that, with probability tending
to one, $\tilde F_n$ belongs to the allowed class, for all large $n$,
and is a consistent estimate of $F_0$.

%
\begin{lemma}
\label{lemmaconsistencyplugin}
Let, under the conditions of Lemma \ref{lemmaexistencesolution},
$F^{(n)}\in B(F_0,\eta)$ be the solution of
\[
\phi\bigl(t;h_1^{(n)},h_2^{(n)},h^{(n)},F
\bigr)=0,\qquad t\in[0,M], %
\]
where $\phi$ is defined by (\ref{fundamentaleq}), and where
$h^{(n)}\in C[S]$, $h_1^{(n)}\in C[S_1]$ an $h_2^{(n)}\in C[S_2]$ are
nonnegative functions which have bounded continuous derivatives on the
supports $S$, $S_1$ and $S_2$, respectively, with finite limits
approaching the boundary, respectively.
Furthermore, let
\[
\bigl\llVert h_j^{(n)}-h_{0j}\bigr\rrVert
\longrightarrow0\quad\mbox{and}\quad\bigl\llVert\bigl(h_j^{(n)}
\bigr)'-h_{0j}'\bigr\rrVert
\longrightarrow0,\qquad j=1,2, %
\]
where, as before, $\llVert \cdot\rrVert $ denotes the supremum norm
on $C[S_j]$.
Finally, let
%
\begin{equation}\label{L1-convergence}
\qquad\bigl\llVert h^{(n)}-h_0\bigr\rrVert
_S\longrightarrow0\quad{and}\quad\bigl\llVert
\partial_j h^{(n)}-\partial_j h_0
\bigr\rrVert_S\longrightarrow0,
\qquad j=1,2,
\end{equation}
where $\llVert \cdot\rrVert _S$ is defined by (\ref{h-norm}).
Then $F^{(n)}\to F_0$ in the supremum metric, as $n\to\infty$, and
$F^{(n)}$ is strictly increasing on $[0,M]$ and satisfies
$F^{(n)}(t)\in[0,1]$, $t\in[0,M]$, for all large $n$.
\end{lemma}

We still need to show that the estimates $\tilde h_{nj}$ of $h_{0j}$
and $\tilde h_n$ of $h_0$ have the properties of $h_j$ and $h$, as
defined in Lemma \ref{lemmaconsistencyplugin}.

%
\begin{lemma}
\label{consistencydensityestimates}
Let the conditions of Theorem \ref{thasympnorsep} be satisfied and
let the estimates $\tilde h_{nj}$ of and $\tilde h_n$ satisfy Condition
\ref{condition2}. Then
\[
\llVert\tilde h_{nj}-h_{0j}\rrVert\stackrel{p}
\longrightarrow0\quad\mbox{and}\quad\bigl\llVert\tilde h_{nj}'-h_{0j}'
\bigr\rrVert\stackrel{p}\longrightarrow0,\qquad j=1,2. %
\]
Moreover,
%
\begin{equation}
\label{L1-convergence-inprob} \qquad\llVert\tilde h_n-h_0\rrVert
_S\stackrel{p}\longrightarrow0\quad{and}\quad\llVert
\partial_j\tilde h_n-\partial_j
h_0\rrVert_S\stackrel{p}\longrightarrow0,\qquad
j=1,2.
\end{equation}
\end{lemma}

We now get the following result.

%
\begin{lemma}
\label{lemmalinearinteq}
Let the conditions of Theorem \ref{thasympnorsep} be satisfied and
let, for small $\eta>0$, $F=\tilde F_n\in B(F_0,\eta)$ be the
solution of the equation
\[
\phi(t;\tilde h_{n1},\tilde h_{n2},\tilde
h_n,F)=0,\qquad t\in[0,M], %
\]
where $\phi$ is defined by (\ref{fundamentaleq}). Moreover, let
$\llVert
\cdot\rrVert $ denote the supremum norm. Then:
\begin{longlist}[(iii)]
\item[(i)] With probability tending to one, $\tilde F_n$ is strictly
increasing on $[0,M]$, and satisfies $\tilde F_n(t)\in[0,1]$, $t\in
[0,M]$, for all large $n$. Hence, with probability tending to one,
$\tilde F_n$ coincides with the MSLE for large $n$ and
\[
\llVert\tilde F_n-F_0\rrVert\stackrel{p}
\longrightarrow0,\qquad n\to\infty. %
\]
\item[(ii)]
\[
\llVert\tilde F_n-F_0\rrVert=O_p
\bigl(n^{-2/5}\sqrt{\log n} \bigr),\qquad n\to\infty. %
\]
\item[(iii)]
\[
\llVert\tilde F_n-\bar F_n\rrVert=O_p
\bigl(n^{-4/5}\log n \bigr), %
\]
where $\bar F_n$ is the solution in $F$ of the linear integral equation
(\ref{linearasympequation1}).
\end{longlist}
\end{lemma}

We did now in principle solve our problem, since we have shown that
$\tilde F_n$ is locally asymptotically equivalent to the solution $\bar
F_n$ of a \textit{linear} integral equation. Since $\tilde F_n$ coincides
with the MSLE for large $n$, the MSLE is also locally asymptotically
equivalent with $\bar F_n$.
However, to get an explicit expression for the bias and variance of the
MSLE, we now study a still simpler ``toy estimator,'' which turns out
also to be locally asymptotically equivalent to the MSLE.

%
\begin{lemma}
\label{lemmatoy}
Let the toy estimator $F=F_n^{\toy}$ be defined as the solution of the equation
%
\begin{eqnarray}\label{toy-equation}
\qquad&&\bigl\{F(t)-F_0(t)\bigr\}\nonumber
\\
&&\quad{}\times \biggl
\{1+d_{F_0}(t) \biggl\{\int_{u<t}\frac
{g(u,t)}{F_0(t)-F_0(u)}
\,dt +\int_{u>v}\frac{g(t,u)}{F_0(u)-F_0(t)} \,du \biggr\} \biggr\}
\nonumber\\[-8pt]\\[-8pt]\nonumber
&&\qquad =\frac{\tilde h_{n1}(t)\{1-F_0(t)\}-\tilde h_{n2}(t)F_0(t)}{\{
1-F_0(t)\}g_1(t)+F_0(t)g_2(t)}
\nonumber
\\
&&\quad\qquad{} +d_{F_0}(t) \biggl\{\int
_{u<t}\frac
{\tilde h_n(u,t)}{F_0(t)-F_0(u)} \,du-\int_{u>v}
\frac{\tilde
h_n(t,u)}{F_0(u)-F_0(t)} \,du \biggr\}.\nonumber
\end{eqnarray}
Then, under the conditions of Theorem \ref{thasympnorsep},
\[
\sqrt{nb_n} \biggl\{F_n^{\toy}(v)-F_0(v)-
\frac{\beta
_1(v)b_n^2}{2\sigma_1(v)} \biggr\}\stackrel{{\mathcal
D}}\longrightarrow N \bigl(0,\sigma
(v)^2 \bigr), %
\]
where $b_n\asymp n^{-1/5}$ and $\beta_1(v)$, $\sigma_1(v)$ and
$\sigma(v)$ are defined as in Theorem \ref{thasympnorsep}.
\end{lemma}

%
\begin{remark}
In Lemma \ref{lemmatoy}, a toy estimator is introduced, which plays a
similar role as the toy estimator in the study of the ordinary MLE for
interval censoring, introduced in \cite{Gr91} and \cite{GrWe92}
(the term ``toy estimator'' was coined by Jon Wellner). It is called a
toy estimator because we cannot use it in an actual sample, since $F_0$
is unknown (and is in fact the object we want to estimate). Actually,
the solution $\bar F_n$ of the linear integral equation (\ref
{linearasympequation1}) in part (iii) of Lemma \ref
{lemmalinearinteq} is also a toy estimator in this sense (but does
not produce explicit expressions for the expectation and variance of
the asymptotic distribution).
\end{remark}

Lemma \ref{lemmatoy2} shows that the solution of the linear integral
equation (\ref{linearasympequation1}) is equivalent in first order
to the toy estimator of Lemma \ref{lemmatoy}, apart from a
deterministic bias term.

%
\begin{lemma}
\label{lemmatoy2}
Let, under the conditions of Theorem \ref{thasympnorsep},
$F_n^{\toy}$ solve equation~(\ref{toy-equation}) of Lemma \ref
{lemmatoy} and let $\bar F_n$ be the solution of the linear integral
equation~(\ref{linearasympequation1}). Then
\begin{eqnarray*}
\bar F_n(t)&=&F_n^{\toy}(t)+d_{F_0}(t)
\biggl\{\int_{u=0}^t\frac{\gamma
_n(u)g(u,t)}{F_0(t)-F_0(u)} \,du+\int
_{u=t}^M\frac{\gamma
_n(u)g(t,u)}{F_0(u)-F_0(t)} \,du \biggr
\}
\\
&&{}+O_p \bigl(n^{-1/2} \bigr),
\end{eqnarray*}
where
\[
\gamma_n(u)=\frac{\beta_1(u)b_n^2}{\sigma_1(u)} %
\]
and $\beta_1(u)$ and $\sigma_1(u)$ are defined as in Theorem \ref
{thasympnorsep}.
\end{lemma}

Theorem \ref{thasympnorsep} now follows from Lemma \ref
{lemmatoy2} and the asymptotic equivalence of~$\bar F_n$ with the MSLE.

\section{Concluding remarks and open problems}
In the preceding, it was shown that, under the so-called separation
hypothesis, the MSLE locally converges to the underlying distribution
function at rate $n^{-2/5}$, if we use bandwidths $b_n\asymp n^{-1/5}$
in the estimates $\tilde h_{nj}$ and $\tilde h_n$. The asymptotic
(normal) distribution was also determined. The results can be used to
construct a two-sample likelihood ratio test, of the same type as the
test, discussed in \cite{piet12} for the current status model, but
this is not done in the present paper. It is also possible to use the
results to construct pointwise bootstrap confidence intervals, as is
done in \cite{pietgeurt14b} and \cite{pietgeurt14} for the
current status model. In that case, it might be advisable to use \textit{undersmoothing}, and work with bandwidths of order $n^{-\alpha}$,
where $1/3<\alpha<1/5$, as is done in \cite{pietgeurt14b} and
\cite{pietgeurt14}. In this way one gets rid of the bias and it is
expected that the SMLE and MSLE will then be very similar, since their
asymptotic variances are the same, which implies that their asymptotic
(normal) limits will also be the same.

If the separation hypothesis does not hold, which means that we can
have arbitrarily small observation intervals, the asymptotic behavior
of the MSLE is still unknown. In this situation the local asymptotic
limit distribution for the ordinary MLE is also still unknown, although
it is conjectured that the rate $n^{-1/3}$, holding under the
separation hypothesis, is improved to the rate $(n\log n)^{-1/3}$ in
this case. There even exists a conjectured limit distribution in this
case, put forward in \cite{Gr91} (see also \cite{pietgeurt14} and
\cite{GrWe92}). Supporting evidence for this conjecture is given in a
simulation study in \cite{piettom11}, but a proof is still missing.
The latter paper also gives simulation results for the SMLE, and the
asymptotic variance of the SMLE is the same as that of the MSLE, but
the asymptotic bias is different, just as in the current status model.
The bias of the SMLE is considerably simpler than the bias of the MSLE.
The asymptotic behavior of the SMLE again has to be deduced from an
associated integral equation; this is further discussed in \cite
{pietgeurt14} and \cite{piettom11}.

It is possible to extend the theory to the situation that there are
more than two observation times than just $T_i$ and $U_i$ or to the
so-called mixed case (see, e.g., \cite{schick00}), where there are a
random number of observation times per unobservable event $X_i$.
However, since this leads to further complications in defining the
integral equations, we did not do this in the present paper.

\begin{appendix}\label{app}
\section*{Appendix}

\begin{pf*}{Proof of Lemma \ref{lemmacharacterization}}
First suppose that the conditions (\ref{fenchel1}) and (\ref
{fenchel2}) are satisfied. Then we cannot have $\hat F_n(t)=0$ for $t$
in an interval where $\tilde h_{n1}(t)>0$ or $\int_{u=0}^t\tilde
h_n(u,t)>0$, since otherwise $\ell(\hat F_n)=-\infty$. Similarly, we
cannot have $\hat F_n(t)=1$ for $t$ in an interval where $\tilde
h_{n2}(t)>0$ or $\int_{u=t}^M\tilde h_n(t,u)>0$.

Since the criterion function $F\mapsto\ell(F)$ is concave in $F$, we get
%
\begin{equation}
\label{concavitycondition} \ell(F)-\ell(\hat F_n)\le\int_0^M
\nabla_{\hat F_n}(u)\bigl\{F(u)-\hat F_n(u)\bigr\} \,du,
\end{equation}
where we use the facts that the integrals defining $\ell(F)$ are all
nonpositive. Note that this is similar to the relation (1.11) in \cite
{GrWe92}.

By (\ref{fenchel2}),
\[
\int_0^M\nabla_{\hat F_n}(u)\hat
F_n(u) \,du=0, %
\]
and hence
\[
\int_0^M\nabla_{\hat F_n}(u)\bigl\{F(u)-
\hat F_n(u)\bigr\} \,du=\int_0^M
\nabla_{\hat F_n}(u) F(u) \,du. %
\]
If $F=1_{[t,\infty)}$ for some $t\in[0,M)$, we get by (\ref{fenchel1}),
\[
\int_0^M\nabla_{\hat F_n}(u) F(u) \,du=
\int_t^M\nabla_{\hat F_n}(u) \,du\le0.
\]
So we also get, for subdistribution functions of the type
\[
F=\sum_{i=1}^k\alpha_i1_{[t_i,\infty)},\qquad
0\le t_1<\cdots, t_k\le M, \alpha_i
\in(0,1), \sum_{i=1}^k\alpha
_i\le1, %
\]
that
\[
\int_0^M\nabla_{\hat F_n}(u) F(u) \,du =
\sum_{i=1}^k\alpha_i\int
_{t_i}^M\nabla_{\hat F_n}(u) \,du\le0.
\]
Since we can approximate any subdistribution function $F$ on $[0,M)$ by
subdistribution functions of this type, this implies
$\ell(F)\le\ell(\hat F_n)$, for all subdistribution functions $F$.

Conversely, suppose that $\hat F_n$ maximizes $\ell(F)$. Then we must
have, if $t\in(0,M)$, $F=1_{[t,\infty)}$
and $\delta\in(0,1)$,
\[
\int_{v=t}^M\nabla_{\hat F_n}(v) \,dv=\lim
_{\delta\downarrow
0}\delta^{-1} \bigl\{\ell\bigl((1-\delta)\hat
F_n+\delta F\bigr)-\ell(\hat F_n) \bigr\}\le0
\]
(using the concavity of $\ell$ for the existence of the limit), and
hence (\ref{fenchel1}) has to be satisfied for $\hat F_n$. Moreover,
defining $F_{\delta}$ by
\[
F_{\delta}(t)=(1+\delta)\hat F_n(t)\wedge1,\qquad t\in[0,M],
\]
we find
\[
\lim_{\delta\to0}\frac{\ell(F_{\delta})-\ell(\hat
F_n)}{\delta}=0,
\]
since the limit has to be nonpositive, if we let $\delta$ tend to
zero, either from above or from below.\vadjust{\goodbreak}

We have
\[
0=\lim_{\delta\to0}\frac{\ell(F_{\delta})-\ell(\hat
F_n)}{\delta}= \int_0^M
\nabla_{\hat F_n}(u)\hat F_n(u) \,du,
\]
so (\ref{fenchel2}) must hold.

Suppose $\hat F_n$ has a jump at $t\in(0,M)$ and suppose $\nabla
_{\hat F_n}(t-)>0$. Define
\[
F_{\delta}(u)=\cases{ \hat F_n(u), &\quad$u<t-\delta$,
\vspace*{3pt}\cr
\hat F_n(t), &\quad$u\in[t-\delta,t)$,
\vspace*{3pt}\cr
\hat F_n(u), &
\quad$u\in[t,M]$.} %
\]
Then
\[
\int\nabla_{\hat F_n}(u)\bigl\{F_{\delta}(u)-\hat F_n(u)
\bigr\} \,du>0,
\]
for small $\delta>0$, a contradiction. Hence, we must have: $\nabla
_{\hat F_n}(t-)\le0$. If $\nabla_{\hat F_n}(t)<0$, we define
\[
F_{\delta}(u)=\cases{ \hat F_n(u), &\quad$u<t$,
\vspace*{3pt}\cr
\hat
F_n(t-), &\quad$u\in[t,t+\delta)$,
\vspace*{3pt}\cr
\hat F_n(u), &
\quad$u\in[t+\delta,M]$,} %
\]
and then again:
\[
\int\nabla_{\hat F_n}(u)\bigl\{F_{\delta}(u)-\hat F_n(u)
\bigr\} \,du>0,
\]
for small $\delta>0$, a contradiction, so we must have: $\nabla
_{\hat F_n}(t)\ge0$, implying we must have
%
\begin{equation}
\label{jumpassumption} \nabla_{\hat F_n}(t-)\le0\le\nabla_{\hat F_n}(t).
\end{equation}

On the other hand, we have by the continuity of $\tilde h_{nj}$,
$j=1,2$ and $\tilde h_n$:
\begin{eqnarray*}
&& \nabla_{\hat F_n}(t)-\nabla_{\hat F_n}(t-)
\\
&&\qquad =\frac{\tilde
h_{n1}(t)}{\hat F_n(t)}-
\frac{\tilde h_{n1}(t)}{\hat
F_n(t-)}-\frac{\tilde h_{n2}(t)}{1-\hat F_n(t)}+\frac{\tilde
h_{n2}(t)}{1-\hat F_n(t-)}
\\
&&\quad\qquad{}+\int_{v=0}^u \frac{\tilde h_{n}(v,t)}{\hat F_n(t)-\hat F_n(v)} \,dv
-\int
_{v=0}^u \frac{\tilde h_{n}(v,t)}{\hat F_n(t-)-\hat F_n(v)} \,dv
\\
&&\quad\qquad{}-\int_{v=t}^M \frac{\tilde h_{n}(t,v)}{\hat F_n(v)-\hat F_n(t)} \,dv
+\int
_{v=t}^M \frac{\tilde h_{n}(t,v)}{\hat F_n(v)-\hat F_n(t-)} \,dv
\\
&&\qquad =\frac{\tilde h_{n1}(t)}{\hat F_n(t)}-\frac{\tilde h_{n1}(t)}{\hat
F_n(t-)}-\frac{\tilde h_{n2}(t)}{1-\hat F_n(t)}+\frac{\tilde
h_{n2}(t)}{1-\hat F_n(t-)}
\\
&&\qquad =\tilde h_{n1}(t)\frac{\hat F_n(t-)-\hat F_n(t)}{\hat F_n(t-)\hat
F_n(t)} +\tilde h_{n2}(t)
\frac{\hat F_n(t-)-\hat F_n(t)}{\{1-\hat F_n(t-)\}\{
1-\hat F_n(t)\}}<0,
\end{eqnarray*}
contradicting (\ref{jumpassumption}). The conclusion is that we must
have: $\hat F_n(t-)=\hat F_n(t)$, for $t\in(0,M)$.

Finally, suppose (\ref{increasepoint}) is satisfied for a point $t\in
(0,M)$ and suppose $\nabla_{\hat F_n}(t)>0$. Then, by the continuity
of $\hat F_n$, there also exists a neighborhood of $t$ such that
$\nabla_{\hat F_n}(u)>0$ for $u$ in this neighborhood. We now define a
perturbation $F_{\delta}$ of $\hat F_n$ by
\[
F_{\delta}(u)=\cases{ \hat F_n(u), &\quad$u<t-\delta$,
\vspace*{3pt}\cr
\hat F_n(t+\delta), &\quad$u\in[t-\delta,t+\delta)$,
\vspace*{3pt}\cr
\hat
F_n(u), &\quad$u\in[t+\delta,M]$.} %
\]
Then we have for sufficiently small $\delta>0$:
\[
\int\nabla_{\hat F_n}(u)\bigl\{F_{\delta}(u)-\hat F_n(u)
\bigr\} \,du>0,
\]
contradicting
\[
\int\nabla_{\hat F_n}(u)\bigl\{F_{\delta}(u)-\hat F_n(u)
\bigr\} \,du\le0.
\]
If (\ref{increasepoint}) is satisfied for a point $t\in(0,M)$ and
$\nabla_{\hat F_n}(t)<0$, we define the perturbation $F_{\delta}$ of
$\hat F_n$ by
\[
F_{\delta}(u)=\cases{ \hat F_n(u), &\quad$u<t-\delta$,
\vspace*{3pt}\cr
\hat F_n(t-\delta), &\quad$u\in[t-\delta,t+\delta)$,
\vspace*{3pt}\cr
\hat
F_n(u), &\quad$u\in[t+\delta,M]$} %
\]
and get a contradiction in the same way. So, if (\ref{increasepoint})
is satisfied for a point $t\in(0,M)$, we must have:
\[
\nabla_{\hat F_n}(t)=0. %
\]
This proves the left-hand side of (\ref{fenchel3}).

Furthermore,
%
\begin{eqnarray}\label{increaseintegral}
&& \int_0^M\nabla_{\hat F_n}(u)
\hat F_n(u) \,du\nonumber
\\
&&\qquad  = \biggl[-\hat F_n(t)\int
_{u=t}^M\nabla_{\hat F_n}(u) \,du
\biggr]_{t=0}^M +\int_{t=0}^M
\int_{u=t}^M\nabla_{\hat F_n}(u) \,du \,d\hat
F_n(t)
\nonumber\\[-8pt]\\[-8pt]\nonumber
&&\qquad =\hat F_n(0)\int_{u=0}^M
\nabla_{\hat F_n}(u) \,du+\int_{t=0}^M\int
_{u=t}^M\nabla_{\hat F_n}(u) \,du \,d\hat
F_n(t)
\nonumber
\\
&&\qquad =\int_{t=0}^M\int_{u=t}^M
\nabla_{\hat F_n}(u) \,du \,d\hat F_n(t),\nonumber
\end{eqnarray}
implying by (\ref{fenchel1}) that
\[
\int_{u=t}^M\nabla_{\hat F_n}(u) \,du=0,
\]
for points $t$ satisfying (\ref{increasepoint}), since otherwise the
right-hand side of (\ref{increaseintegral}) would be strictly negative.
\end{pf*}

\begin{pf*}{Proof of Theorem \ref{thconsistency}}
By the\vspace*{2pt} assumption on the kernel estimates and the observation density
$g$, we may assume that $\tilde h_n$ satisfies (\ref
{separationconditionestimate}), for some $\varepsilon>0$, and all
large $n$. Let the function $\psi$ be defined by
%
\begin{eqnarray}\label{criterionfunction-psi}
\psi(F;h_1,h_2,h) &=& \int
h_1(t)\log F(t) \,dt+\int h_2(t)\log\bigl\{1-F(t)\bigr\}
\,dt
\nonumber\\[-8pt]\\[-8pt]\nonumber
&&{} +\int h(t,u)\log\bigl\{F(u)-F(t)\bigr\} \,dt \,du.
\end{eqnarray}
Then we must have, if $h_j=\tilde h_{nj}$, $j=1,2$ and $h=\tilde h_n$,
\[
\lim_{\varepsilon\downarrow0}\varepsilon^{-1} \bigl\{\psi\bigl((1-
\varepsilon)\hat F_n+\varepsilon F_0;h_1,h_2,h
\bigr)- \psi(\hat F_n;h_1,h_2,h) \bigr\}
\le0.
\]
This implies (see also (4.20) in \cite{GrWe92}):
%
\begin{eqnarray}\label{consistency-ineq}
&& \int\frac{F_0(t)}{\hat F_n(t)} \tilde
h_{n1}(t) \,dt+\int
\frac
{1-F_0(t)}{1-\hat F_n(t)} \tilde h_{n2}(t) \,dt
\nonumber\\[-8pt]\\[-8pt]\nonumber
&&\qquad{} +\int\frac
{F_0(u)-F_0(t)}{\hat F_n(u)-\hat F_n(t)} \tilde
h_n(t,u) \,dt \,du\le1.
\end{eqnarray}
Fix a small $\delta\in[0,M/2]$ and let the intervals $A_{\delta}$
and $B_{\delta}$ be defined by
\[
A_{\delta}=[\delta,M],\qquad B_{\delta}=[0,M-\delta]. %
\]
Then, arguing as in \cite{GrWe92}, Part II, Chapter~4, we find that
there exists an $M>0$ such that for all $n$,
\[
\sup_{t\in A_{\delta}}1\bigm/\hat F_n(t)+\sup_{t\in B_{\delta}}1\bigm/
\bigl\{ 1-\hat F_n(t)\bigr\}\le M. %
\]
We also cannot have that $\hat F_n(u_{k_n})-\hat F_n(t_{k_n})\to0$,
for a sequence of points $(t_{k_n},u_{k_n})\in C_{\delta}$. For
suppose, if necessary by taking a subsequence, that $t_{k_n}\to t_0$
and $u_{k_n}\to u_0$. The $u_0-t_0\ge\varepsilon+\delta$. By the
vague convergence of $\hat F_n$ to $F$, there are continuity points
$t_1$ and $u_1$ such that
$t_0<t_1<u_1<u_0$, $u_1-t_1\ge\frac12\delta+\varepsilon$, and
$\hat F_n(t_1)\to F(t_1)$ and $\hat F_n(u_1)\to F(u_1)$. Moreover, since
\[
\hat F_n(u_1)-\hat F_n(t_1)\le
\hat F_n(u_{k_n})-\hat F_n(t_{k_n}),
\]
for large $n$, we must have: $F(u_1)-F(t_1)=0$. We then would get that
there exists a rectangle
$[t_1,t_2]\times[u_2,u_1]$ such that $u_2-t_2>\varepsilon$ and
\begin{eqnarray*}
&&\liminf_{n\to\infty}\int_{[t_1,t_2]\times[u_2,u_1]}
\frac
{F_0(u)-F_0(t)}{\hat F_n(u)-\hat F_n(t)} \tilde h_n(t,u) \,dt \,du
\\
&&\qquad\ge K\bigl\{F_0(u_2)-F_0(t_2)
\bigr\}\int_{[t_1,t_2]\times[u_2,u_1]}h_0(t,u) \,dt \,du,
\end{eqnarray*}
for any $K>0$, contradicting (\ref{consistency-ineq}). So, we may also
assume that
\[
\inf_{(t,u)\in C_{\delta}}\bigl\{\hat F_n(u)-\hat
F_n(t)\bigr\}\ge\frac{1}{M}, %
\]
for all $n$.

As in \cite{GrWe92}, Part II, Chapter~4, we have by the Helly
compactness theorem that the exists a set of probability one, such that
for each $\omega$ in this set the sequence $(\hat F_n(\cdot;\omega
))$ has a subsequence $\hat F_{n_k}(\cdot;\omega)$ which converges
vaguely to a subdistribution function $F=F(\cdot;\omega)$.
By the vague convergence of $\hat F_{n_k}(\cdot;\omega)$ to $F$, we
now get
\begin{eqnarray*}
&&\int_{A_{\delta}} \frac{F_0(t)}{\hat F_{n_k}(t;\omega)} \tilde h_{n1}(t)
\,dt+\int_{B_{\delta}}\frac{1-F_0(t)}{1-\hat
F_{n_k}(t;\omega)} \tilde h_{n2}(t) \,dt
\\
&&\quad{} +\int_{C_{\delta}} \frac{F_0(u)-F_0(t)}{\hat F_{n_k}(u;\omega
)-\hat F_{n_k}(t;\omega)} \tilde h_n(t,u) \,dt
\,du
\\
&&\qquad \to\int_{A_{\delta}}\frac{F_0(t)}{F(t)} h_{01}(t) \,dt+\int
_{B_{\delta}}\frac{1-F_0(t)}{1-F(t)} h_{02}(t) \,dt
\\
&&\quad\qquad{} +\int
_{C_{\delta}} \frac{F_0(u)-F_0(t)}{F(u)-F(t)} h_0(t,u) \,dt \,du,\qquad n\to
\infty.
\end{eqnarray*}
By monotone convergence, we now also have
%
\begin{eqnarray}\label{limit-boundedness}
&&\int_{[0,M]}\frac{F_0(t)}{F(t)}
h_{01}(t) \,dt+\int_{[0,M]}\frac
{1-F_0(t)}{1-F(t)}
h_{02}(t) \,dt\nonumber
\\
&&\quad{}  +\int_{[0,M]^2}\frac{F_0(u)-F_0(t)}{F(u)-F(t)}
h_0(t,u) \,dt \,du
\nonumber\\[-8pt]\\[-8pt]\nonumber
&&\qquad =\lim_{\delta\downarrow0} \biggl\{\int_{A_{\delta}}
\frac{F_0(t)}{F(t)} h_{01}(t) \,dt+\int_{B_{\delta}}
\frac{1-F_0(t)}{1-F(t)} h_{02}(t) \,dt
\\
&&\hspace*{95pt}{} +\int_{C_{\delta}}
\frac{F_0(u)-F_0(t)}{F(u)-F(t)} h_0(t,u) \,dt \,du \biggr\}\le1.\nonumber
\end{eqnarray}
Suppose $F(t)\ne F_0(t)$ for some $t\in[0,M/2]$. Then there exist a
$u\in(t,M)$ such that $h_0(t,u)>0$ and
\[
\frac{F_0(t)^2}{F(t)}+\frac{\{1-F_0(u)\}^2}{1-F(u)}+\frac{\{
F_0(u)-F_0(t)\}^2}{F(u)-F(t)}>1,
\]
since
\begin{eqnarray*}
\frac{F_0(t)^2}{x}+\frac{\{1-F_0(u)\}^2}{1-y}+\frac{\{F_0(u)-F_0(t)\}
^2}{y-x} \cases{ =1, &
\quad$F_0(t)=x$, $F_0(u)=y$,
\vspace*{3pt}\cr
>1, &\quad otherwise}
\end{eqnarray*}
(see also (4.27) in \cite{GrWe92}). By the continuity of $F_0$ and
the monotonicity and right continuity of $F$ there exist therefore also
$h>0$ such that
\[
\frac{F_0(t')^2}{F(t')}+\frac{\{1-F_0(u')\}^2}{1-F(u')}+\frac{\{
F_0(u')-F_0(t')\}^2}{F(u')-F(t')}>1, %
\]
if $t'\in[t,t+h]$ and $u'\in[u,u+h]$. This implies
\begin{eqnarray*}
&&\int_{[0,M]}\frac{F_0(t)}{F(t)} h_{01}(t) \,dt+\int
_{[0,M]}\frac
{1-F_0(t)}{1-F(t)} h_{02}(t) \,dt
\\
&&\quad{}+\int
_{[0,M]^2}\frac{F_0(u)-F_0(t)}{F(u)-F(t)} h_0(t,u) \,dt \,du
\\
&&\qquad =\int_{[0,M]}\frac{F_0(t)^2}{F(t)}g_1(t) \,dt+\int
_{[0,M]}\frac{\{
1-F_0(t)\}^2}{1-F(t)} g_2(t) \,dt
\\
&&\quad\qquad{} +\int
_{[0,M]^2}\frac{\{F_0(u)-F_0(t)\}^2}{F(u)-F(t)} g(t,u) \,dt \,du
\\
&&\qquad =\int_{[0,M]^2} \biggl\{\frac{F_0(t)^2}{F(t)}+\frac{\{1-F_0(t)\}^2}{1-F(t)}
+\frac{\{F_0(u)-F_0(t)\}^2}{F(u)-F(t)} \biggr\} g(t,u) \,dt \,du>1,
\end{eqnarray*}
in contradiction with (\ref{limit-boundedness}). So, we must have
$F(t)=F_0(t)$ if $t\in[0,M/2]$. A~similar argument yields
$F(t)=F_0(t)$ if $t\in[M/2,M)$.

So, for each $\omega$ outside a set of probability zero, the sequence
$(\hat F_n(\cdot;\omega))$ has a subsequence which converges weakly
to $F_0$. This implies that $\hat F_n(t)$ converges almost surely to
$F_0(t)$ for each $t\in[0,M)$.
The uniformity of the convergence on subintervals follows from the
continuity of $F_0$.
\end{pf*}

\begin{pf*}{Proof of Lemma \ref{lemmaexistencesolution}}
We will use the line of argument of the proof of the implicit function
theorem 10.2.1 in \cite{dieudonne1969}. We define the function $\bar
\phi$
by
%
\begin{equation}
\label{bar-phi} \bigl[\bar\phi(h_1,h_2,h,F) \bigr](t)=
\phi(t,h_1,h_2,h,F),\qquad t\in[0,M],
\end{equation}
so $\bar\phi$ maps $E$ to $C[0,M]$.
The derivative of $\bar\phi$ w.r.t. $F$, is given by the function
\begin{eqnarray*}
&& \bigl[ \bigl[\partial_4\bar\phi(h_1,h_2,h,F)
\bigr](A) \bigr](t)
\\
&&\qquad\stackrel{\mathrm{def}}=- \bigl\{h_1(t)+h_2(t)
\bigr\}A(t)
\\
&&\quad\qquad{}+\bigl\{ 1-2F(t)\bigr\} \biggl\{\int_{v=0}^t
\frac{h(v,t)}{F(t)-F(v)} \,dv -\int_{u=t}^M
\frac{h(t,u)}{F(u)-F(t)} \,du \biggr\}A(t)
\\
&&\quad\qquad{}-F(t)\bigl\{1-F(t)\bigr\}
\\
&&\qquad\qquad{}\times \biggl\{\int_{u=0}^t
\frac{h(u,t) \{
A(t)-A(u) \}}{\{F(t)-F(u)\}^2} \,du+\int_{u=t}^M
\frac{h(t,u)\{
A(t)-A(u)\}}{\{F(u)-F(t)\}^2} \,du \biggr\},
\end{eqnarray*}
where $A\in C[0,M]$. Note that the right-hand side is well defined for
$F\in B(F_0,\eta)$ and small $\eta>0$, since $h(t,u)=0$ if
$u-t<\varepsilon$, and since $F_0$ has a nonvanishing derivative on
$[0,M]$, implying that $F_0(u)-F_0(t)$ stays away from zero if $u-t\ge
\varepsilon$.

We now define the open set $U=U_{\delta}$ of functions $(h_1,h_2,h)$ by
%
\begin{eqnarray}\label{U-delta}
U_{\delta} &=& \bigl\{(h_1,h_2,h)\in
C[S_1]\times C[S_2]\times C[S]\dvtx
\nonumber\\[-8pt]\\[-8pt]\nonumber
&&\hspace*{5pt} \max\bigl\{\llVert
h_1-h_{01}\rrVert,\llVert h_2-h_{02}
\rrVert,\llVert h-h_0\rrVert_S \bigr\}<\delta\bigr\}.
\end{eqnarray}
There exists a $\delta>0$ such that for $(h_1,h_2,h)\in U_{\delta}$,
\begin{eqnarray*}
&& \bigl\llVert\bar\phi(h_1,h_2,h,F_1)-
\bar\phi(h_1,h_2,h,F_2)- \bigl[
\partial_4\bar\phi(h_{01},h_{02},h_0,F_0)
\bigr](F_1-F_2)\bigr\rrVert
\\
&&\qquad \le\varepsilon '\llVert F_1-F_2\rrVert,
\end{eqnarray*}
if $F_1,F_2\in B(F_0,\eta)$, where $\varepsilon'>0$ can be made
arbitrarily small by making $\delta$ small, using the definition of
differentiability in Banach spaces.

The equation
\[
\bigl[\partial_4\bar\phi(F_0;h_{01},h_{02},h_0,F_0)
\bigr](A)=0 %
\]
only has the trivial solution $A\equiv0$ in $C[0,M]$. This is seen in
the following way. Suppose there exists a solution in $A\in C[0,M]$
such that $ [\partial_4\bar\phi
(F_0;h_{01},h_{02},\break  h_0, F_0) ](A)=0$ and $A(s)>0$ for some $s\in
[0,M]$. Then also $\max_{s\in[0,M]}A(s)>0$. Suppose the maximum is
attained at $t\in[0,M]$. Then
\begin{eqnarray*}
&& \bigl[ \bigl[\partial_4\bar\phi(h_{01},h_{02},h_0,F_0)
\bigr](A) \bigr](t)
\\
&&\qquad =- \bigl\{h_{01}(t)+h_{02}(t) \bigr\}A(t)
\\
&&\quad\qquad{} +\bigl
\{1-2F_0(t)\bigr\} \biggl\{\int_{v=0}^t
\frac{h_0(v,t)}{F_0(t)-F_0(v)} \,dv -\int_{u=t}^M
\frac{h_0(t,u)}{F_0(u)-F_0(t)} \,du \biggr\}A(t)
\\
&&\quad\qquad{}-F_0(t)\bigl\{1-F_0(t)\bigr\}
\\
&&\qquad\qquad{}\times  \biggl\{\int
_{u=0}^t\frac{h_0(u,t) \{
A(t)-A(u) \}}{\{F_0(t)-F_0(u)\}^2} \,du+\int
_{u=t}^M\frac
{h_0(t,u)\{A(t)-A(u)\}}{\{F_0(u)-F_0(t)\}^2} \,du \biggr\}
\\
&&\qquad =- \bigl\{g_1(t)\bigl\{1-F_0(t)\bigr
\}+g_2(t)F_0(t) \bigr\}A(t)
\\
&&\quad\qquad{}-F_0(t)\bigl\{1-F_0(t)\bigr\}
\\
&&\qquad\qquad{}\times \biggl\{\int
_{u=0}^t\frac{h_0(u,t) \{
A(t)-A(u) \}}{\{F_0(t)-F_0(u)\}^2} \,du+\int
_{u=t}^M\frac
{h_0(t,u)\{A(t)-A(u)\}}{\{F_0(u)-F_0(t)\}^2} \,du \biggr\}
\\
&&\qquad \le- \bigl\{g_1(t)\bigl\{1-F_0(t)\bigr
\}+g_2(t)F_0(t) \bigr\}A(t)<0,
\end{eqnarray*}
using $g_1(t)\{1-F_0(t)\}+g_2(t)F_0(t)>0$, in contradiction with the assumption
\[
\bigl[\partial_4\bar\phi(h_{01},h_{02},h_0,F_0)
\bigr](A)=0. %
\]
We similarly get a contradiction if we assume that $A(t)<0$ for some
$t\in[0,M]$ (similar arguments were used for the integral equation,
studied in \cite{GeGr96}). This shows that $\partial_4\bar\phi
(h_{01},h_{02},h_0,F)$ is a linear homeomorphism of $C[0,M]$ onto
$C[0,M]$ and that we can in fact use arguments of the type used in the
proof of the implicit function theorem in Banach spaces, as given, for
example, in \cite{dieudonne1969}, Theorem 10.2.1.

Denoting (as in the proof of Theorem 10.2.1 of \cite{dieudonne1969})
the linear mapping $\partial_4\bar\phi(h_{01},h_{02},h_0,F_0)$ by
$T_0$ and its inverse by $T_0^{-1}$, we find that
\begin{eqnarray*}
\bigl\llVert T_0^{-1}\cdot\bigl\{\bar\phi
(h_1,h_2,h,F_1)-\bar\phi
(h_1,h_2,h,F_2) \bigr\}-
(F_1-F_2 )\bigr\rrVert
&\le&\varepsilon'
\bigl\llVert T_0^{-1}\bigr\rrVert\llVert
F_1-F_2\rrVert
\\
&\le&\tfrac12\llVert F_1-F_2
\rrVert,
\end{eqnarray*}
so we have a \textit{contraction}, and this implies that the equation
\[
F=F-T_0^{-1}\cdot\,\bar\phi(\cdot,F;h_1,h_2,h)
\]
has a unique solution $F\in B(F_0,\eta)$ which can be obtained by
successive approximations, if we take the balls around $h_{0j}$ and
$h_0$, to which $h_j$ and $h$ belong, respectively, sufficiently small,
using a result like 10.1.1 in \cite{dieudonne1969}. This, in turn,
implies that the equation
\[
\bar\phi(h_1,h_2,h,F)=0 %
\]
has a unique solution in $F\in B(F_0,\eta)$, for $(h_1,h_2,h)\in
U_{\delta}$ and small $\delta$.
\end{pf*}

\begin{pf*}{Proof of Lemma \ref{lemmaderivativesolution}}
If $\phi(\cdot;h_1,h_2,h,F)=0$, we have
\begin{eqnarray*}
&& \bigl\{1-F(t)\bigr\}h_1(t)-F(t)h_2(t)
\\
&&\quad {}
+F(t)\bigl\{1-F(t)\bigr\} \biggl\{\int_{u\dvtx (u,t)\in S}
\frac
{h(u,t)}{F(t)-F(u)} \,du -\int_{u\dvtx (t,u)\in S}\frac{h(t,u)}{F(u)-F(t)} \,du
\biggr
\}
\\
&&\qquad =0.
\end{eqnarray*}
Note that the differentiability proerties of $h$, $h_1$ and $h_2$ and
the fact that $F$ solves the integral equation imply that we can
differentiate $F$ too.
Differentiation w.r.t. $t$, and defining $f=F'$, yields:
%
\begin{eqnarray}\label{derivativeeq}
&&\bigl\{1-F(t)\bigr\}h_1'(t)-F(t)h_2'(t)-f(t)
\bigl\{h_1(t)+h_2(t)\bigr\}\nonumber
\\
&&\quad{} +\bigl\{1-2F(t)\bigr\}\nonumber
\\
&&\qquad{}\times  f(t)\biggl\{\int_{u\dvtx (u,t)\in S}
\frac{h(u,t)}{F(t)-F(u)} \,du -\int_{u\dvtx (t,u)\in S}\frac{h(t,u)}{F(u)-F(t)} \,du
\biggr
\}\nonumber
\\
&&\quad{}+F(t)\bigl\{1-F(t)\bigr\}
\nonumber\\[-8pt]\\[-8pt]\nonumber
&&\qquad{}\times  \biggl\{\int_{u\dvtx (u,t)\in S}
\frac{\partial
_2h(u,t)}{F(t)-F(u)} \,du -\int_{u\dvtx (t,u)\in S}\frac{\partial
_1h(t,u)}{F(u)-F(t)} \,du
\biggr\}
\\
&&\quad{}-f(t)F(t)\bigl\{1-F(t)\bigr\} \nonumber
\\
&&\qquad{}\times \biggl\{\int_{u\dvtx (u,t)\in S}
\frac{h(u,t)}{\{
F(t)-F(u)\}^2} \,du +\int_{u\dvtx (t,u)\in S}\frac{h(t,u)}{\{F(u)-F(t)\}^2} \,du
\biggr\}\nonumber
\\
&&\qquad =0.\nonumber
\end{eqnarray}
Temporarily replacing $F$ by $F_0$, and $h_j$ and $h$ by $h_{0j}$ and
$h_0$, respectively, we would obtain
\begin{eqnarray*}
&&\bigl\{1-F_0(t)\bigr\}h_{01}'(t)-F_0(t)h_{02}'(t)-f(t)
\bigl\{\bigl\{1-F_0(t)\bigr\} g_1(t)+F_0(t)g_2(t)
\bigr\}
\\
&&\quad{}+F_0(t)\bigl\{1-F_0(t)\bigr\}
\\
&&\qquad{}\times \biggl\{\int
_{u\dvtx (u,t)\in S}\frac{\partial
_2h_0(u,t)}{F_0(t)-F_0(u)} \,du -\int_{u\dvtx (t,u)\in S}
\frac{\partial_1h_0(t,u)}{F_0(u)-F_0(t)} \,du \biggr\}
\nonumber
\\
&&\quad{}-f(t)F_0(t)\bigl\{1-F_0(t)\bigr\}
\\
&&\qquad{}\times \biggl\{\int
_{u\dvtx (u,t)\in S}\frac
{h_0(u,t)}{\{F_0(t)-F_0(u)\}^2} \,du +\int_{u\dvtx (t,u)\in S}
\frac{h_0(t,u)}{\{F_0(u)-F_0(t)\}^2} \,du \biggr\}
\\
&&\qquad =0,
\end{eqnarray*}
that is,
\begin{eqnarray*}
&&f_0(t) \biggl\{\bigl\{\bigl\{1-F_0(t)\bigr
\}g_1(t)1_{S_1}(t)+F_0(t)g_2(t)1_{S_2}(t)
\bigr\}
\\
&&\hspace*{17pt}\quad{}+F_0(t)\bigl\{1-F_0(t)\bigr\}
\biggl\{\int_{u\dvtx (u,t)\in
S}\frac{h_0(u,t)}{\{F_0(t)-F_0(u)\}^2} \,du
\\
&&\hspace*{120pt}{}  +\int
_{u\dvtx (t,u)\in S}\frac{h_0(t,u)}{\{F_0(u)-F_0(t)\}^2} \,du \biggr\} \biggr\}
\\
&&\qquad =\bigl\{1-F_0(t)\bigr\}h_{01}'(t)1_{S_1}(t)-F_0(t)h_{02}'(t)1_{S_2}(t)
\\
&&\quad\qquad{}+F_0(t)\bigl\{1-F_0(t)\bigr\}
\\
&&\qquad\qquad{}\times \biggl\{
\int_{u\dvtx (u,t)\in S}\frac
{\partial_2h_0(u,t)}{F_0(t)-F_0(u)} \,du -\int_{u\dvtx (t,u)\in S}
\frac{\partial_1h_0(t,u)}{F_0(u)-F_0(t)} \,du \biggr\}.
\end{eqnarray*}
This means that the coefficient of $f(t)$ in equation (\ref
{derivativeeq}) stays away from zero if $F$ belongs to
a sufficiently small ball $B(F_0,\eta)$ around $F_0$, and
$(h_1,h_2,h)\in U_{\delta}$ for small $\delta>0$. Denoting this
coefficient by $c(t,h_1,h_2,h,F)$, we get the equation
%
\begin{eqnarray}\label{uniformboundednessder}
f(t)&=&\frac{1}{c(t,h_1,h_2,h,F)}\nonumber
\\
&&{}\times  \biggl\{
\bigl\{1-F(t)\bigr\}
h_1'(t)-F(t)h_2'(t)
\nonumber\\[-8pt]\\[-8pt]\nonumber
&&\hspace*{18pt}{} +F(t)\bigl\{1-F(t)\bigr\}
\\
&&\hspace*{28pt}{}\times \biggl\{\int
_{u=0}^t\frac{\partial_2h(u,t)}{F(t)-F(u)} \,du -\int
_{u=t}^M\frac{\partial_1h(t,u)}{F(u)-F(t)} \,du \biggr\} \biggr\}.\nonumber
\end{eqnarray}
The statement of the lemma now follows.
\end{pf*}

\begin{pf*}{Proof of Lemma \ref{lemmaconsistencyplugin}}
By Lemma \ref{lemmaexistencesolution}, we have that $F^{(n)}$ tends
to $F_0$ in the supremum norm on $C[0,M]$, since for any $\eta$ we can
choose a $\delta>0$ such that
\[
\bigl\llVert F^{(n)}-F_0\bigr\rrVert<\eta, %
\]
if $\llVert h_j^{(n)}-h_{0j}\rrVert <\delta$ and $\llVert
h^{(n)}-h_0\rrVert <\delta$.

Using Lemma \ref{lemmaderivativesolution}, we get that, under the
conditions of the lemma, that $F^{(n)}$ is differentiable with a
bounded derivative $f^{(n)}$; see (\ref{uniformboundednessder}).
Specifically, (\ref{uniformboundednessder}) yields
\begin{eqnarray*}
f^{(n)}(t)&\stackrel{\mathrm{def}}=&\bigl(F^{(n)}
\bigr)'(t)
\\
&=& \biggl\{F^{(n)}(t)\bigl\{1-F^{(n)}(t)\bigr\}
\\
&&\hspace*{3pt}{}\times \biggl\{\int
_{u=0}^t\frac
{\partial_2h^{(n)}(u,t)}{F^{(n)}(t)-F^{(n)}(u)} \,du -\int
_{u=t}^M\frac{\partial_1h^{(n)}(t,u)}{F^{(n)}(u)-F^{(n)}(t)} \,du \biggr\}
\\
&&\hspace*{77pt}{}+\bigl\{1-F^{(n)}(t)\bigr\}
\bigl(h_1^{(n)} \bigr)'(t)-F^{(n)}(t)
\bigl(h_2^{(n)} \bigr)'(t) \biggr\}
\\
&&{}\Bigm/c
\bigl(t,h_1^{(n)},h_2^{(n)},h^{(n)},F^{(n)}
\bigr),
\end{eqnarray*}
where $c(t,h_1^{(n)},h_2^{(n)},h^{(n)},F^{(n)})$ stays away from zero,
as $n\to\infty$.
The corresponding density $f_0$ of the underlying model similarly has
the representation
%
\begin{eqnarray}\label{f0-relation}
\qquad f_0(t)&=& \biggl\{F_0(t)\bigl\{1-F_0(t)\bigr\}\nonumber
\\
&&\hspace*{3pt}{}\times \biggl\{\int_{u\dvtx (u,t)\in S}
\frac
{\partial_2h_0(u,t)}{F_0(t)-F_0(u)} \,du -\int_{u\dvtx (t,u)\in S}\frac
{\partial_1h_0(t,u)}{F_0(u)-F_0(t)} \,du
\biggr\}
\nonumber
\\
&&\hspace*{131pt}{}+\bigl
\{1-F_0(t)\bigr\}h_{01}'(t)-F_0(t)h_{02}'(t)
\biggr\}
\\
&&{}\Bigm/c_0(t),\nonumber
\end{eqnarray}
where $c_0(t)$ is given by
\[
c_0(t)=g_1(t)\bigl\{1-F_0(t)\bigr
\}+g_2(t)F_0(t). %
\]

By
\[
\bigl\llVert h^{(n)}_j-h_{0j}\bigr\rrVert\to0,
\qquad\bigl\llVert h^{(n)}-h_0\bigr\rrVert\to0,\qquad\bigl
\llVert F^{(n)}-F_0\bigr\rrVert\to0, %
\]
and (\ref{L1-convergence}), we now get
\[
\sup_{t\in[0,M]}\bigl\llvert c \bigl(t,h_1^{(n)},h_2^{(n)},h^{(n)},F^{(n)}
\bigr)-c_0(t)\bigr\rrvert\to0. %
\]
Again using (\ref{L1-convergence}), we also get
\[
\bigl\llVert f^{(n)}-f_0\bigr\rrVert\to0, %
\]
that is, $f^{(n)}$ converges to $f_0$ in the supremum norm. Since $f_0$
stays away from zero on $[0,M]$, this means that $F^{(n)}$ is strictly
increasing on $[0,M]$ for all sufficiently large $n$.

Furthermore, since $\phi(t,h^{(n)}_1,h^{(n)}_2,h^{(n)},F^{(n)})=0$,
we get for large $n$, and $t$ in a right neighborhood of $0$,
\begin{eqnarray*}
\hspace*{-2pt}&& F^{(n)}(t)
\\
\hspace*{-2pt}&&\quad  =\frac{h_1^{(n)}(t)}{h_1^{(n)}(t)+h_2^{(n)}(t)+\{
1-F^{(n)}(t)\}\int_{u=t}^M (\fracc{h^{(n)}(t,u)}{F^{(n)}(u)-F^{(n)}(t)})
\,du}
\\
\hspace*{-2pt}&&\quad \ge0,
\end{eqnarray*}
since, by the convergence of $F^{(n)}$ to $F_0$, we may assume
$1-F^{(n)}(t)>0$ for $t$ in a neighborhood of $0$, and since
$h_j^{(n)}$ and $h^{(n)}$ are nonnegative.\vspace*{2pt}

Likewise, if $F_0(M)=1$, we have, for $t$ in a small left neighborhood
of $M$,
\begin{eqnarray*}
&& F^{(n)}(t)\bigl\{1-F^{(n)}(t)\bigr\}\int_{v=0}^t
\frac
{h^{(n)}(v,t)}{F^{(n)}(t)-F^{(n)}(v)} \,dv
\\
&&\qquad = \bigl\{h_1^{(n)}(t)+h_2^{(n)}(t)
\bigr\}F^{(n)}(t)-h_1^{(n)}(t)
\\
&&\quad\qquad{} +F^{(n)}(t)
\bigl\{1-F^{(n)}(t)\bigr\}\int_{v=t}^M
\frac
{h^{(n)}(v,t)}{F^{(n)}(t)-F^{(n)}(v)} \,dv
\\
&&\qquad =h_2^{(n)}(t)F^{(n)}(t)-h_1^{(n)}(t)
\bigl\{1-F^{(n)}(t)\bigr\},
\end{eqnarray*}
and hence, for $t$ in a small left neighborhood of $M$,
\[
\bigl\{1-F^{(n)}(t)\bigr\} \biggl\{F^{(n)}(t)\int
_{v=0}^t\frac
{h^{(n)}(v,t)}{F^{(n)}(t)-F^{(n)}(v)} \,dv+h_1^{(n)}(t)
\biggr\} =h_2^{(n)}(t)F^{(n)}(t),
\]
implying that, for all large $n$, $1-F^{(n)}(t)\ge0$ for $t$ in a
neighborhood of $M$. This will a fortiori hold if $F_0(M)<1$.
This shows that, for all large $n$ and all $t\in[0,M]$, $F^{(n)}(t)\in[0,1]$.
\end{pf*}

\begin{pf*}{Proof of Lemma \ref{consistencydensityestimates}}
Since we use boundary kernels near the boundary of $[0,M]$, $\tilde
h_{nj}(t)$ is a consistent estimate of $h_{0j}(t)$ for each $t\in S_j$.
For if $t\in[b_n,M-b_n]\cap S_j$ we just have
\[
{\mathbb E}\tilde h_{n1}(t)={\mathbb E}\Delta_{11}
K_{b_n}(t-T_1)=\int K_{b_n}(t-u)h_{01}(u)
\,du=h_{01}(t)+O\bigl(n^{-2/5}\bigr), %
\]
where the remainder term is uniform in $t\in[b_n,M-b_n]\cap S_j$.
Since $b_n\downarrow0$, we have $b_n<\varepsilon$ for all large $n$,
where $\varepsilon$ is the ``separation parameter'' of Condition \ref
{condition1}, and hence the boundary kernels are only relevant for
$\tilde h_{n1}$ in a neighborhood of $0$ and for $\tilde h_{n2}$ in a
neighborhood of $M$.

If $t\in[0,b_n]$, we have
\begin{eqnarray*}
{\mathbb E}\tilde h_{n1}(t)&=&\alpha(t/b_n)E\Delta
_{11} K_{b_n}(t-T_1)+\beta(t/b_n){
\mathbb E} \biggl\{\Delta_{11} \frac
{t-T_1}{b_n}K_{b_n}(t-T_1)
\biggr\}
\\
&=&\alpha(t/b_n)\int_{u=0}^M
K_{b_n}(t-u)h_{01}(u) \,du
\\
&&{}+\beta(t/b_n)\int
_{u=0}^M \frac{t-u}{b_n}K_{b_n}(t-u)h_{01}(u)
\,du
\\
&=&h_{01}(t)\int_{u=-1}^{t/b_n} \bigl\{\alpha
(t/b_n)K(u)+\beta(t/b_n)K(u) \bigr\} \,du+O
\bigl(n^{-2/5}\bigr)
\\
&=&h_{01}(t)+O\bigl(n^{-2/5}\bigr),
\end{eqnarray*}
again uniformly for $t\in[0,b_n]$. A similar computation can be made
for $\tilde h_{n2}$ if $t\in[M-b_n,M]$. Since
\[
\sup_{t\in S_1}\bigl\llvert\tilde h_{n1}(t)-{\mathbb E}
\tilde h_{n1}(t)\bigr\rrvert=O_p \bigl(n^{-2/5}
\sqrt{\log n} \bigr), %
\]
we now get the uniform convergence in probability of $\tilde h_{n1}$ to
$h_{01}$ on $S_1$, and similarly we have uniform convergence in
probability of $\tilde h_{n2}$ to $h_{02}$ on $S_2$.\vspace*{2pt}

Next, we consider the derivative of $\tilde h_{n1}(t)$. If $t\in
[b_n,M-b_n]\cap S_1$, we just have
\begin{eqnarray*}
{\mathbb E}\tilde h_{n1}'(t)&=&\frac{d}{dt}{
\mathbb E}\Delta_{11} K_{b_n}'(t-T_1)
\\
&=&
\int\frac{d}{dt}K_{b_n}(t-u)h_{01}(u) \,du
\\
&=&b_n^{-1}\int K'(u)h_{01}(t-b_nu)
\,du
\\
&=&b_n^{-1}\int K'(u) \biggl\{h_{01}(t)-b_nuh_{01}'(t)
+\frac12b_n^2 u^2h_{01}'(t)-
\frac16 b_n^3 u^3h_{01}''(t)
\biggr\} \,du
\\
&&{} +o \bigl(b_n^2 \bigr)
\\
&=&h_{01}'(t)+O (b_n )=h_{01}'(t)+O
\bigl(n^{-1/5} \bigr),
\end{eqnarray*}
again uniformly in $t$. Since
\[
\sup_{t\in[b_n,M-b_n]\cap S_1}\bigl\llvert\tilde h_{n1}'(t)-{
\mathbb E}\tilde h_{n1}'(t)\bigr\rrvert=O_p
\bigl(n^{-1/5}\sqrt{\log n} \bigr), %
\]
we only have to consider what happens near the boundary.

In treating the boundary kernels, we denote for simplicity $b_n$ by
$b$. If $t\in[0,b]$, we have
\begin{eqnarray*}
{\mathbb E}\tilde h_{n1}(t)&=&\alpha\biggl(\frac{t}{b} \biggr)
\int_{x=0}^{t+b} K_b(t-x)h_{01}(x)
\,dx
\\
&&{} +\beta\biggl(\frac{t}b \biggr)\int_{x=0}^{t+b}
\frac{t-x}b K_b(t-x)h_{01}(x) \,dx.
\end{eqnarray*}
This can also be written
\begin{eqnarray*}
{\mathbb E}\tilde h_{n1}(t)&=&\int_{-1}^{t/b}
\biggl\{\alpha\biggl(\frac{t}{b} \biggr)K(u)+\beta\biggl(
\frac{t}{b} \biggr)uK(u) \biggr\}h_{01}(t-bu) \,du.
\end{eqnarray*}
We write this in the form
\begin{eqnarray*}
&& {\mathbb E}\tilde h_{n1}(t)
\\
&&\qquad =h_{01}(t)
\\
&&\quad\qquad{} +\int
_{-1}^{t/b} \biggl\{\alpha\biggl(\frac{t}{b}
\biggr)K(u)+\beta\biggl(\frac{t}{b} \biggr)uK(u) \biggr\}\int
_{t-bu}^{t}(w-t+bu)h_{01}''(w)
\,dw \,du,
\end{eqnarray*}
using a second-order Taylor development of $h_{01}$ with the integral
remainder term. Hence,
\begin{eqnarray*}
\frac{d}{dt}{\mathbb E}\tilde h_{n1}(t)
&=&h_{01}'(t)+
\frac{1}{b} \biggl\{ \alpha\biggl(\frac{t}{b} \biggr)K \biggl(
\frac{t}{b} \biggr)+\beta\biggl(\frac{t}{b} \biggr)
\frac{t}{b}K \biggl(\frac{t}{b} \biggr) \biggr\}\int
_{0}^{t}wh_{01}''(w)
\,dw
\\
&&{}+\frac{1}{b}\int_{-1}^{t/b}
\biggl\{\alpha' \biggl(\frac
{t}{b} \biggr)K(u)+\beta
' \biggl(\frac{t}{b} \biggr)uK(u) \biggr\}
\\
&&\quad{}\times  \int
_{t-bu}^{t}(w-t+bu)h_{01}''(w)
\,dw \,du
\\
&&{}-\int_{-1}^{t/b} \biggl\{\alpha
\biggl(\frac
{t}{b} \biggr)K(u)+\beta\biggl(\frac{t}{b}
\biggr)uK(u) \biggr\}\int_{t-bu}^{t}h_{01}''(w)
\,dw \,du
\\
&&{}+\int_{-1}^{t/b} \biggl\{\alpha
\biggl(\frac
{t}{b} \biggr)K(u)+\beta\biggl(\frac{t}{b}
\biggr)uK(u) \biggr\}bu h_{01}''(t) \,du
\\
&=& h_{01}'(t)+O(b),\qquad b\downarrow0.
\end{eqnarray*}
Note that, by Condition \ref{condition2}, the functions $\alpha$,
$\beta$, $\alpha'$ and $\beta'$ are bounded on $[0,1]$.
We also have
\begin{eqnarray*}
\frac{d}{dt}{\mathbb E}\tilde h_{n1}(t)&=&\frac{1}{b}
\biggl\{\alpha' \biggl(\frac{t}{b} \biggr)\int
_{x=0}^{t+b} K_b(t-x)h_{01}(x)
\,dx
\\
&&\hspace*{12pt}{}+\beta' \biggl(\frac{t}b \biggr)\int
_{x=0}^{t+b} \frac{t-x}b K_b(t-x)h_{01}(x)
\,dx \biggr\}
\\
&&{} +\alpha\biggl(\frac{t}{b} \biggr)\int_{x=0}^{t+b}
\frac
{d}{dt}K_b(t-x)h_{01}(x) \,dx
\\
&&{}+\beta\biggl(
\frac{t}b \biggr)\int_{x=0}^{t+b}
\frac{t-x}b \frac{d}{dt}K_b'(t-x)h_{01}(x)
\,dx
\\
&&{}+\frac{1}{b}\beta\biggl(\frac{t}b \biggr)\int
_{x=0}^{t+b} K_b(t-x)h_{01}(x)
\,dx,
\end{eqnarray*}
so
\[
{\mathbb E}\tilde h_{n1}'(t)=\frac{d}{dt}{\mathbb
E}\tilde h_{n1}(t)=h_{01}'(t)+O(b)=h_{01}'(t)+O
\bigl(n^{-1/5} \bigr).
\]
Since we have
\[
\sup_{t\in[0,b_n]}\bigl\llvert\tilde h_{n1}'(t)-{
\mathbb E}\tilde h_{n1}'(t)\bigr\rrvert=O_p
\bigl(n^{-1/5} \bigr), %
\]
we now also get that
\[
\sup_{u\in S_1}\bigl\llvert\tilde h_{n1}'(t)-h_{01}'(t)
\bigr\rrvert=o_p(1). %
\]
The other cases can be treated in a similar way.
\end{pf*}

\begin{pf*}{Proof of Lemma \ref{lemmalinearinteq}}
Part (i) is an immediate consequence of Lemma \ref
{lemmaconsistencyplugin}.

(ii)~We get, again using the approach of the implicit function theorem
10.2.1 in Banach spaces of \cite{dieudonne1969}, denoting the
derivative w.r.t. $(h_1,h_2,h)$ by $D_1$ and the derivative w.r.t. $F$
by $D_2$:
%
\begin{eqnarray}
\label{differentiationinteq}
\llVert\tilde F_n-F_0\rrVert&=&\bigl
\llVert\bigl[D_2\bar\phi(h_{01},h_{02},h_0,F_0)^{-1}
\circ D_1\bar\phi(h_{01},h_{02},h_0,F_0)
\bigr]\nonumber
\\
&&\hspace*{74pt}{}\times (\tilde h_{n1}-h_{01},\tilde h_{n2}-h_{02},
\tilde h_n-h_0)\bigr\rrVert
\\
&&{} +o_p \bigl(\bigl
\llVert(\tilde h_{n1}-h_{01},\tilde h_{n2}-h_{02},
\tilde h_n-h_0)\bigr\rrVert\bigr),\nonumber
\end{eqnarray}
where the norm $\llVert \cdot\rrVert $ on the left-hand side and the
first norm on
the right-hand side denote the supremum norm on $C[0,M]$ and the norm
in the $o_p$-term denotes the norm
\[
\bigl\llVert(h_1,h_2,h)\bigr\rrVert=\max\bigl\{
\llVert h_1\rrVert,\llVert h_2\rrVert,\llVert h\rrVert
_S \bigr\}, %
\]
where the first two norms denote again the supremum norm and the third
norm $\llVert \cdot\rrVert _S$ is defined by
(\ref{h-norm}).

By well-known results in density estimation, we have, if $b_n\asymp n^{-1/5}$,
\[
\max\bigl(\llVert\tilde h_{n1}-h_{01}\rrVert,\llVert
\tilde h_{n2}-h_{02}\rrVert\bigr)=O_p
\bigl(n^{-2/5}\sqrt{\log n} \bigr). %
\]
The boundary kernels ensure that the rates are not spoiled by what
happens at the boundary.
So, we have to determine the rate of convergence of $\llVert \tilde
h_n-h_0\rrVert
_S$. We get
\begin{eqnarray*}
&& \int_{u\dvtx (u,t)\in S}\bigl\llvert\tilde h_n(u,t)-h_0(u,t)
\bigr\rrvert \,du
\\
&&\qquad \le M^{1/2} \biggl\{\int_{u\dvtx (u,t)\in S} \bigl\{
\tilde h_n(u,t)-h_0(u,t) \bigr\}^2 \,du \biggr
\}^{1/2}
\end{eqnarray*}
and
\[
\biggl\{\int_{u\dvtx (u,t)\in S} \bigl\{\tilde h_n(u,t)-E
\tilde h_n(u,t) \bigr\}^2 \,du \biggr\}^{1/2}=O_p
\bigl(n^{-2/5}\sqrt{\log n} \bigr),
\]
uniformly in $t$.
For the bias we get, if $b_n<u<u+\varepsilon\le t<M-b_n$,
\begin{eqnarray*}
&& {\mathbb E}\tilde h_n(u,t)-h_0(u,t)
\\
&&\qquad ={\mathbb
E}K_{b_n}(u-T_1)K_{b_n}(t-U_1)\Delta
_{12}-h_0(u,t)
\\
&&\qquad =\int K_{b_n}\bigl(u-t'\bigr)K_{b_n}
\bigl(t-u'\bigr)h_0\bigl(t',u'
\bigr) \,dt' \,du'-h_0(u,t)
\\
&&\qquad =\int K(v)K(w)h_0(u-b_nw,t-b_nv) \,dv
\,dw-h_0(u,t)
\\
&&\qquad =O \bigl(b_n^2 \bigr).
\end{eqnarray*}
The use of the boundary kernels ensures that the bias is also of order
$O(b_n^2)$ is $u<b_n$ or $t>M-b_n$.
The conclusion is
%
\begin{equation}
\label{parameterupperbound} \bigl\llVert(\tilde h_{n1}-h_{01},\tilde
h_{n2}-h_{02},\tilde h_n-h_0)\bigr
\rrVert=O_p \bigl(n^{-2/5}\sqrt{\log n} \bigr).
\end{equation}

The derivative $D_2$ was computed in the proof of Lemma \ref
{lemmaexistencesolution} (denoted by $\partial_4$ there) and the
derivative $D_1$ is given by
\begin{eqnarray*}
&& \bigl[ \bigl[D_1\bar\phi(h_{01},h_{02},h_0,F_0)
\bigr](A) \bigr](t)
\\
&&\qquad =B_1(t)\bigl\{1-F_0(t)\bigr\}-B_2(t)F_0(t)
\nonumber
\\
&&\quad\qquad{} +F_0(t)\bigl\{1-F_0(t)\bigr\}
\biggl\{\int_{v=0}^t \frac
{B(v,t)}{F_0(t)-F_0(v)} \,dv -\int
_{u=t}^M \frac{B(t,u)}{F_0(u)-F_0(t)} \,du \biggr\},
\end{eqnarray*}
where $B_1$, $B_2$ and $B$ are of the form
\begin{eqnarray*}
B_1&=& h_1-h_{01},\qquad B_2=h_2-h_{02},
\\
B&=&h-h_0. %
\end{eqnarray*}

Hence, defining $\bar F_n$ by
%
\begin{eqnarray}\label{linearsolution}
\bar F_n-F_0 &=&- \bigl[D_2\bar
\phi(h_{01},h_{02},h_0,F_0)^{-1}
\circ D_1\bar\phi(h_{01},h_{02},h_0,F_0)
\bigr]
\nonumber\\[-8pt]\\[-8pt]\nonumber
&&\hspace*{6pt}{}\times (\tilde h_{n1}-h_{01},\tilde h_{n2}-h_{02},
\tilde h_n-h_0),
\end{eqnarray}
we get that $F=\bar F_n$ is the solution of the linear integral equation
\begin{eqnarray*}
&& \bigl[D_2\bar\phi(h_{01},h_{02},h_0,F_0)
\bigr] (F-F_0 )
\\
&&\qquad =- \bigl[D_1\bar\phi
(h_{01},h_{02},h_0,F_0) \bigr](
\tilde h_{n1}-h_{01},\tilde h_{n2}-h_{02},
\tilde h_n-h_0), %
\end{eqnarray*}
which, letting $A=F-F_0$ and $(B_1,B_2,B)=(\tilde h_{n1}-h_{01},\tilde
h_{n2}-h_{02},\tilde h_n-h_0)$, boils down to the equation
%
\begin{eqnarray}\label{lininteq2}
 && \bigl\{h_{01}(t)+h_{02}(t) \bigr\}A(t)\nonumber
\\
&&\quad{}-
\bigl\{1-2F_0(t)\bigr\} \biggl\{\int_{v=0}^t
\frac{h_0(v,t)}{F_0(t)-F_0(v)} \,dv -\int_{u=t}^M
\frac{h_0(t,u)}{F_0(u)-F_0(t)} \,du \biggr\} A(t)
\nonumber
\\
&&\quad{} +F_0(t)\bigl\{1-F_0(t)\bigr\} \nonumber
\\
&&\qquad{}\times \biggl\{\int
_{u=0}^t\frac{h_0(u,t) \{
A(t)-A(u) \}}{\{F_0(t)-F_0(u)\}^2} \,du+\int
_{u=t}^M\frac
{h_0(t,u)\{A(t)-A(u)\}}{\{F_0(u)-F_0(t)\}^2} \,du \biggr\}\hspace*{-20pt}
\\
&&\qquad =B_1(t)\bigl\{1-F_0(t)\bigr\}-B_2(t)F_0(t)
\nonumber
\\
&&\qquad\quad{}+F_0(t)\bigl\{1-F_0(t)\bigr\}\nonumber
\\
&&\qquad\qquad{}\times \biggl\{\int_{v=0}^t \frac
{B(v,t)}{F_0(t)-F_0(v)} \,dv -\int
_{u=t}^M \frac{B(t,u)}{F_0(u)-F_0(t)} \,du \biggr\}.\nonumber
\end{eqnarray}
We have
\begin{eqnarray*}
&& \bigl\{h_{01}(t)+h_{02}(t) \bigr\}A(t)
\\
&&\quad{} -\bigl
\{1-2F_0(t)\bigr\} \biggl\{\int_{v=0}^t
\frac{h_0(v,t)}{F_0(t)-F_0(v)} \,dv -\int_{u=t}^M
\frac{h_0(t,u)}{F_0(u)-F_0(t)} \,du \biggr\}A(t)
\\
&&\qquad = \bigl\{g_1(t)F_0(t) (t)+g_2(t)\bigl
\{1-F_0(t)\bigr\} \bigr\}A(t)
\\
&&\quad\qquad{}-\bigl\{1-2F_0(t)\bigr\}
\bigl\{g_2(t)-g_1(t)\bigr\}A(t)
\\
&&\qquad = \bigl\{\bigl\{1-F_0(t)\bigr\}g_1(t)+F_0(t)g_2(t)
\bigr\}A(t).
\end{eqnarray*}
Furthermore,
\begin{eqnarray*}
&&\int_{u=0}^t\frac{h_0(u,t) \{A(t)-A(u) \}}{\{
F_0(t)-F_0(u)\}^2} \,du+\int
_{u=t}^M\frac{h_0(t,u)\{A(t)-A(u)\}}{\{
F_0(u)-F_0(t)\}^2} \,du
\\
&&\qquad =\int_{u=0}^t\frac{g(u,t) \{A(t)-A(u) \}}{F_0(t)-F_0(u)} \,du+\int
_{u=t}^M\frac{g(t,u)\{A(t)-A(u)\}}{F_0(u)-F_0(t)} \,du.
\end{eqnarray*}
Finally,
\begin{eqnarray*}
\hspace*{-3pt}&& h_{01}(t)\bigl\{1-F_0(t)\bigr\}-h_{02}(t)F_0(t)
\\
\hspace*{-3pt}&&\qquad{}+F_0(t)\bigl\{1-F_0(t)\bigr\}
\biggl\{\int_{v=0}^t \frac
{h_0(v,t)}{F_0(t)-F_0(v)} \,dv -\int
_{u=t}^M \frac{h_0(t,u)}{F_0(u)-F_0(t)} \,du \biggr\}
=0.
\end{eqnarray*}
So, we obtain the linear integral equation (\ref
{linearasympequation1}) by dividing both sides of (\ref{lininteq2})
by $\{1-F_0(t)\}g_1(t)+F_0(t)g_2(t)$.\vspace*{1.5pt}

Hence, $\bar F_n$ is the solution of the linear integral equation (\ref
{linearasympequation1}), and by (\ref{differentiationinteq}) and~(\ref{parameterupperbound}),
\[
\llVert\tilde F_n-F_0\rrVert=\llVert\bar
F_n-F_0\rrVert+o_p \bigl(n^{-2/5}
\sqrt{\log n} \bigr). %
\]
Note that we have
\[
\llVert\bar F_n-F_0\rrVert=O_p
\bigl(n^{-2/5}\sqrt{\log n} \bigr) %
\]
by the fact that $D_1\bar\phi(h_{01},h_{02},h_0,F_0)$ and $D_2\bar
\phi(h_{01},h_{02},h_0,F_0)^{-1}$ are bounded linear mappings, and hence
\[
\llVert\bar F_n-F_0\rrVert=O_p \bigl(
\bigl\llVert(\tilde h_{n1}-h_{01},\tilde h_{n2}-h_{02},
\tilde h_n-h_0)\bigr\rrVert\bigr)=O_p
\bigl(n^{-2/5}\sqrt{\log n} \bigr). %
\]
(iii) The function $\tilde F_n$ satisfies the equation $\bar\phi
(\tilde h_{n1},\tilde h_{n2},\tilde h_n,\tilde F_n)=0$, where $\bar
\phi$ is defined by (\ref{bar-phi}). Hence,
\begin{eqnarray*}
&&\tilde h_{n1}(t)\bigl\{1-\tilde F_n(t)\bigr\}-\tilde
h_{n2}(t)\tilde F_n(t)
\\
&&\qquad {} +F(t)\bigl\{1-F(t)\bigr\} \biggl\{\int_{v=0}^t
\frac{\tilde
h_n(v,t)}{\tilde F_n(t)-\tilde F_n(v)} \,dv -\int_{u=t}^M
\frac{\tilde h_n(t,u)}{\tilde F_n(u)-\tilde F_n(t)} \,du \biggr\}
=0.
\end{eqnarray*}
By (ii), we have $\llVert \tilde F_n-F_0\rrVert =O_p(n^{-2/5}\sqrt
{\log n})$, and hence
\begin{eqnarray*}
&&\int_{v=0}^t \frac{\tilde h_n(v,t)}{\tilde F_n(t)-\tilde F_n(v)} \,dv
\\
&&\qquad =\int_{v=0}^t \frac{\tilde h_n(v,t)}{F_0(t)-F_0(v)} \,dv
\\
&&\quad\qquad{}-\int
_{v=0}^t\frac{g(u,t)\{\tilde F_n(t)-F_0(t)-\tilde F_n(u)+F_0(u)\}
}{F_0(t)-F_0(u)} \,du+O_p
\bigl(n^{-4/5}\log n \bigr)
\end{eqnarray*}
and similarly
\begin{eqnarray*}
&&\int_{u=t}^M \frac{\tilde h_n(t,u)}{\tilde F_n(u)-\tilde F_n(t)} \,du
\\
&&\qquad =\int_{u=t}^M \frac{\tilde h_n(t,u)}{F_0(u)-F_0(t)} \,du
\\
&&\quad\qquad{}-\int
_{u=t}^M\frac{g(t,u)\{\tilde F_n(u)-F_0(u)-\tilde F_n(t)+F_0(t)\}
}{F_0(u)-F_0(t)} \,du+O_p
\bigl(n^{-4/5}\log n \bigr).
\end{eqnarray*}
Hence, we get
\begin{eqnarray*}
&&\tilde F_n(t)-F_0(t) +d_{F_0}(t) \biggl\{\int
_{u=0}^t\frac{g(u,t)\{\tilde
F_n(t)-F_0(t)-\tilde F_n(u)+F_0(u)\}}{F_0(t)-F_0(u)} \,du
\\
&&\hspace*{105pt}{} -\int
_{u=t}^M\frac{g(t,u)\{\tilde F_n(u)-F_0(u)-\tilde F_n(t)+F_0(t)\}
}{F_0(u)-F_0(t)} \,du \biggr\}
\nonumber
\\
&&\qquad =\frac{\tilde h_{n1}(t)\{1-F_0(t)\}-\tilde h_{n2}(t)F_0(t)}{\{
1-F_0(t)\}g_1(t)+F_0(t)g_2(t)}
\nonumber
\\
&&\quad\qquad{} +d_{F_0}(t) \biggl\{\int_{u=0}^t
\frac{\tilde
h_n(u,t)}{F_0(t)-F_0(u)} \,du-\int_{u=t}^M
\frac{\tilde
h_n(t,u)}{F_0(u)-F_0(t)} \,du \biggr\}
\\
&&\quad\qquad{}+O_p \bigl(n^{-4/5}\log n
\bigr),
\end{eqnarray*}
uniformly for $t\in[0,M]$, implying
\begin{eqnarray*}
\tilde F_n &=&- \bigl[D_2\bar\phi(h_{01},h_{02},h_0,F_0)^{-1}
\circ D_1\bar\phi(h_{01},h_{02},h_0,F_0)
\bigr]
\\
&&\hspace*{9pt}{}\times (\tilde h_{n1}-h_{01},\tilde h_{n2}-h_{02},
\tilde h_n-h_0)
\\
&&{}+O_p \bigl(n^{-4/5}\log n \bigr)
\\
&=&\bar F_n+O_p \bigl(n^{-4/5}\log n \bigr).
\end{eqnarray*}\upqed
\end{pf*}

\begin{pf*}{Proof of Lemma \ref{lemmatoy}}
We have
\[
\tilde h_{n1}(v)=\frac{1}n\sum_{i=1}^n
K_{b_n}(v-T_i)\Delta_{i1}
\]
and hence
\begin{eqnarray*}
\operatorname{var} \bigl(\tilde h_{n1}(v) \bigr)&=&\frac{1}n\operatorname{var}
\bigl(K_{b_n}(v-T_1)\Delta_{11} \bigr)
\\
&=&\frac{1}n EK_{b_n}(v-T_1)^2 \bigl(
\Delta_{11}-F_0(T_1) \bigr)^2
\\
&\sim&
\frac{F_0(v)\{1-F_0(v)\}g_1(v)}{nb_n}\int K(u)^2 \,du.
\end{eqnarray*}
Likewise,
\[
\operatorname{var} \bigl(\tilde h_{n2}(v) \bigr)\sim\frac{F_0(v)\{
1-F_0(v)\}g_2(v)}{nb_n}\int
K(u)^2 \,du.
\]
Furthermore,
\begin{eqnarray*}
&& \operatorname{covar} \bigl(\tilde h_{n1}(v),\tilde h_{n2}(v)
\bigr)
\\
&&\qquad =\frac{1}n {\mathbb E}K_{b_n}(v-T_1)K_{b_n}(v-U_1)
\bigl(\Delta_{11}-F_0(T_1) \bigr) \bigl(
\Delta_{13}-F_0(U_1) \bigr)
\\
&&\qquad =-\frac{1}n {\mathbb E}K_{b_n}(v-T_1)K_{b_n}(v-U_1)F_0(T_1)F_0(U_1)
\\
&&\qquad
\sim-\frac{F_0(v)^2g(v,v)}{nb_n}\int K(u)^2 \,du=0,
\end{eqnarray*}
using the ``separation condition'' $g(v,v)=0$. So, we obtain
\begin{eqnarray*}
&& \operatorname{var} \biggl(\frac{\{1-F_0(v)\}\tilde h_{n1}(v)-F_0(v)\tilde
h_{n2}(v)}{g_1(v)\{1-F_0(v)\}+F_0(v)g_2(v)} \biggr)
\\
&&\qquad \sim\frac{F_0(v)\{
1-F_0(v)\} (\{1-F_0(v)\}^2g_1(v)+F_0(v)^2g_2(v) )}{
nb_n \{g_1(v)\{1-F_0(v)\}+F_0(v)g_2(v) \}^2}.
\end{eqnarray*}
Furthermore,
\[
\int_{t<v}\frac{\tilde h_n(t,v)}{F_0(v)-F_0(t)} \,dt= n^{-1}\sum
_{i=1}^nK_{b_n}(v-U_i)
\Delta_{i2}\int_{t<v}\frac{
K_{b_n}(t-T_i)}{F_0(v)-F_0(t)} \,dt
\]
and hence
\begin{eqnarray*}
&& \operatorname{var} \biggl(\int_{t<v}\frac{\tilde h_n(t,v)}{F_0(v)-F_0(t)} \,dt
\biggr)
\\
&&\qquad =n^{-1}\operatorname{var} \biggl(K_{b_n}(v-U_1)\Delta
_{12}\int_{t<v}\frac{
K_{b_n}(t-T_1)}{F_0(v)-F_0(t)} \,dt \biggr)
\\
&&\qquad =n^{-1}{\mathbb E}K_{b_n}(v-U_1)^2
\bigl(\Delta_{12}-F_0(U_1)+F_0(T_1)
\bigr)^2 \biggl\{\int_{t<v}\frac{ K_{b_n}(t-T_1)}{F_0(v)-F_0(t)} \,dt
\biggr)^2
\\
&&\qquad =n^{-1}{\mathbb E}K_{b_n}(v-U_1)^2
\bigl\{F_0(U_1)-F_0(T_1) \bigr
\} \bigl\{1-F_0(U_1)+F_0(T_1)
\bigr\}
\\
&&\quad\qquad{}\times  \biggl\{\int_{t<v}\frac{ K_{b_n}(t-T_1)}{F_0(v)-F_0(t)} \,dt \biggr)^2
\\
&&\qquad \sim \frac{1}{nb_n}\int_{t<v} \frac{1-F_0(v)+F_0(t)}{F_0(v)-F_0(t)}
g(t,v) \,dt\int K(u)^2 \,du.
\end{eqnarray*}
Likewise,
\begin{eqnarray*}
&& \operatorname{var} \biggl(\int_{w>v}\frac{\tilde h_n(v,w)}{F_0(w)-F_0(v)} \,dw
\biggr)
\\
&&\qquad \sim\frac{1}{nb_n}\int_{w>v} \frac{1-F_0(w)+F_0(v)}{F_0(w)-F_0(v)} g(v,w)
\,dw\int K(u)^2 \,du.
\end{eqnarray*}
Finally,
\begin{eqnarray*}
&& \operatorname{covar} \biggl(\tilde h_{n2}(v),\int_{t<v}
\frac{\tilde
h_n(t,v)}{F_0(v)-F_0(t)} \,dt \biggr)
\\
&&\qquad =\frac{1}n {\mathbb E}K_{b_n}(v-U_1)^2
\bigl(\Delta_{13}-\bigl\{ 1-F_0(U_1)\bigr\}
\bigr) \bigl(\Delta_{12}-F_0(U_1)+F_0(T_1)
\bigr)
\\
&&\quad\qquad{}\times \int_{t<v}\frac{ K_{b_n}(t-T_1)}{F_0(v)-F_0(t)} \,dt
\\
&&\qquad =-\frac{1}n {\mathbb E}K_{b_n}(v-U_1)^2
\bigl\{1-F_0(U_1)\bigr\} \bigl\{ F_0(U_1)-F_0(T_1)
\bigr\}
\\
&&\hspace*{6pt}\quad\qquad{}\times \int_{t<v}\frac{ K_{b_n}(t-T_1)}{F_0(v)-F_0(t)} \,dt
\\
&&\qquad \sim-\frac{\{1-F_0(v)\}g_2(v)}{nb_n}\int K(u)^2 \,du
\end{eqnarray*}
and similarly
\[
\operatorname{covar} \biggl(\tilde h_{n1}(v),\int_{w>v}
\frac{\tilde
h_n(v,w)}{F_0(w)-F_0(v)} \,dw \biggr) \sim-\frac{F_0(v)g_1(v)}{nb_n}\int
K(u)^2 \,du.
\]

Combining these facts, we obtain that the variance of the right-hand
side of~(\ref{toy-equation}) is given by
\begin{eqnarray*}
\hspace*{-2pt}&&\frac{F_0(v)\{1-F_0(v)\} (\{1-F_0(v)\}
^2g_1(v)+F_0(v)^2g_2(v) )}{
nb_n \{g_1(v)\{1-F_0(v)\}
+F_0(v)g_2(v) \}^2}\int K(u)^2 \,du
\\
\hspace*{-2pt}&&\quad{} +\frac{d_{F_0}(v)^2}{nb_n}
\\
\hspace*{-2pt}&&\qquad{}\times \biggl\{\int_{t<v}
\frac
{g(t,v)}{F_0(v)-F_0(t)} \,dt +\int_{w>v} \frac{g(v,w)}{F_0(w)-F_0(v)}
\,dw-g_1(v)-g_2(v) \biggr\}
\\
\hspace*{-2pt}&&\qquad{}\times \int K(u)^2 \,du
+\frac{2d_{F_0}(v)^2}{nb_n} \bigl\{g_1(v)+g_2(v) \bigr
\}\int K(u)^2 \,du
\\
\hspace*{-2pt}&&\qquad =\frac{F_0(v)\{1-F_0(v)\} (\{1-F_0(v)\}
^2g_1(v)+F_0(v)^2g_2(v) )}{
nb_n \{g_1(v)\{1-F_0(v)\}+F_0(v)g_2(v) \}^2}\int K(u)^2 \,du
\\
\hspace*{-2pt}&&\quad\qquad{} +\frac{d_{F_0}(v)^2}{nb_n} \biggl\{\int_{t<v}
\frac
{g(t,v)}{F_0(v)-F_0(t)} \,dt
\\
\hspace*{-2pt}&&\hspace*{87pt}{} +\int_{w>v} \frac{g(v,w)}{F_0(w)-F_0(v)}
\,dw+g_1(v)+g_2(v) \biggr\} \int K(u)^2 \,du
\\
\hspace*{-2pt}&&\qquad =\frac{d_{F_0}(v)}{
nb_n} \biggl\{1+d_{F_0}(v) \biggl\{\int
_{t<v}\frac
{g(t,v)}{F_0(v)-F_0(t)} \,dt +\int_{w>v}
\frac{g(v,w)}{F_0(w)-F_0(v)} \,dw \biggr\} \biggr\}
\\
\hspace*{-2pt}&&\quad\qquad{}\times \int K(u)^2 \,du.
\end{eqnarray*}
Hence, we get that the asymptotic variance at a fixed interior point
$v$ of the solution $F$ of the equation (\ref{toy-equation}) is given by
\[
\frac{d_{F_0}(v)}{nb_n\sigma_1}\int K(u)^2 \,du,
\]
where $\sigma_1$ is defined by (\ref{sigma1}).

We still have to compute the bias.
We have
\begin{eqnarray*}
&&{\mathbb E}\frac{\{1-F_0(v)\} \{\tilde h_1(v)-h_1(v) \}
-F_0(v) \{\tilde h_2(v)-h_2(v) \}}{g_1(v)\{1-F_0(v)\}
+F_0(v)g_2(v)}
\\
&&\quad{} +d_{F_0}(v){\mathbb E}\int_{t<v}
\frac{\tilde
h(t,v)-h(t,v)}{F_0(v)-F_0(t)} \,dt-d_{F_0}(v){\mathbb E}\int_{u>v}
\frac
{\tilde h_n(v,u)-h(v,u)}{F_0(u)-F_0(v)} \,du
\\
&&\qquad =\frac{\{1-F_0(v)\}h_1''(v)-F_0(v)h_2''(v)}{2 \{g_1(v)\{
1-F_0(v)\}+F_0(v)g_2(v) \}}b_n^2\int u^2K(u) \,du
\\
&&\quad\qquad{} +\frac12b_n^2\,d_{F_0}(v) \biggl\{\int
_{t<v}\frac{(\fraca{\partial^2}{\partial v^2})h(t,v)}{F_0(v)-F_0(t)} \,dt-\int_{u>v}
\frac
{(\fraca{\partial^2}{\partial v^2})h(v,u)}{F_0(u)-F_0(v)} \,du \biggr\}
\\
&&\qquad\qquad{}\times \int
u^2K(u) \,du+o
\bigl(b_n^2 \bigr),
\end{eqnarray*}
where
%
\begin{eqnarray}\label{h1-tilde}
\tilde h_1(t)&=&\int K_{b_n}(t-u)F_0(u)
g_1(u) \,du,
\nonumber\\[-9pt]\\[-9pt]\nonumber
\tilde h_2(t)&=&\int K_{b_n}(t-u)
\bigl\{1-F_0(u)\bigr\} g_2(u) \,du
\end{eqnarray}
and
%
\begin{equation}
\label{h-tilde} \qquad\tilde h(t,u)=\int K_{b_n}(t-v)K_{b_n}(u-w)
\bigl\{F_0(w)-F_0(v)\bigr\} g(v,w) \,dv \,dw.
\end{equation}
Moreover,
\begin{eqnarray*}
&&\bigl\{1-F_0(v)\bigr\}h_1(v)-F_0(v)h_2(v)
\\[-2pt]
&&\quad{} +F_0(v)\bigl\{1-F_0(v)\bigr\} \biggl\{
\int_{t<v}\frac
{h(t,v)}{F_0(v)-F_0(t)} \,dt -\int_{u>v}
\frac{h(v,u)}{F_0(u)-F_0(v)} \,du \biggr\}
\\[-2pt]
&&\qquad =F_0(v)\bigl\{1-F_0(v)\bigr\} \bigl
\{g_1(v)-g_2(v) \bigr\}-F_0(v)\bigl
\{1-F_0(v)\bigr\} \bigl\{g_1(v)-g_2(v) \bigr
\}
\\[-2pt]
&&\qquad =0.
\end{eqnarray*}
This yields the result of the lemma, since we can use the central limit
theorem for i.i.d. random variables on the right-hand side of (\ref
{toy-equation}), using
%
\begin{eqnarray}
\label{sumrepresentation} &&\frac{\tilde h_{n1}(t)\{1-F_0(t)\}-\tilde
h_{n2}(t)F_0(t)}{\{
1-F_0(t)\}g_1(t)+F_0(t)g_2(t)}
\nonumber
\\
&&\quad{}+d_{F_0}(t) \biggl\{\int
_{u<t}\frac
{\tilde h_n(u,t)}{F_0(t)-F_0(u)} \,du-\int_{u>t}
\frac{\tilde
h_n(t,u)}{F_0(u)-F_0(t)} \,du \biggr\}
\nonumber
\\
&&\qquad =n^{-1}\sum_{i=1}^n \biggl\{
\frac{\{1-F_0(t)\} K_{b_n}(t-T_i)\Delta
_{i1}-F_0(t)K_{b_n}(t-U_i)\Delta_{i3}}{\{1-F_0(t)\}
g_1(t)+F_0(t)g_2(t)}
\\
&&\hspace*{38pt}\quad\qquad{} +d_{F_0}(t) K_{b_n}(t-U_i)
\Delta_{i2}\int_{u<t}\frac{ K_{b_n}(u-T_i)}{F_0(t)-F_0(u)} \,du
\nonumber
\\
&&\hspace*{56pt}\quad\qquad{} -d_{F_0}(t)
K_{b_n}(t-T_i)\Delta_{i2}\int
_{u>t}\frac{
K_{b_n}(u-U_i)}{F_0(u)-F_0(t)} \,du \biggr\}.\nonumber
\end{eqnarray}\upqed
\end{pf*}

\begin{pf*}{Proof of Lemma \ref{lemmatoy2}}
The difference between the equation defining the toy estimator and the
solution of the integral equation (\ref{linearasympequation1})
resides in the term
\[
-d_{F_0}(t) \biggl\{\int_{u=0}^t
\frac{g(u,t)\{F(u)-F_0(u)\}
}{F_0(t)-F_0(u)} \,du+\int_{u=t}^M
\frac{g(t,u)\{F(u)-F_0(u)\}
}{F_0(u)-F_0(t)} \,du \biggr\};
\]
see (\ref{linearasympequation1}) and (\ref
{linearasympequation2}). Take $F=F_n^{\toy}$, and consider the
integral within the brackets. An easy computation yields
\[
\int_{u=0}^t\frac{ \{F_n^{\toy}(u)-{\mathbb E}F_n^{\toy}(u)
\}g(u,t)}{F_0(t)-F_0(u)}
\,du=O_p \bigl(n^{-1/2} \bigr).
\]
So, we get
\begin{eqnarray*}
&&d_{F_0}(t) \biggl\{\int_{u=0}^t
\frac{g(u,t)\{F_n^{\toy}(u)-F_0(u)\}
}{F_0(t)-F_0(u)} \,du+\int_{u=t}^M
\frac{g(t,u)\{F_n^{\toy}(u)-F_0(u)\}
}{F_0(u)-F_0(t)} \,du \biggr\}
\\
&&\qquad =d_{F_0}(t) \biggl\{\int_{u=0}^t
\frac{g(u,t)\{{\mathbb
E}F_n^{\toy}(u)-F_0(u)\}}{F_0(t)-F_0(u)} \,du
\\
&&\hspace*{64pt}{}+\int_{u=t}^M
\frac{g(t,u)\{
{\mathbb E}F_n^{\toy}(u)-F_0(u)\}}{F_0(u)-F_0(t)} \,du \biggr\} +O_p \bigl
(n^{-1/2} \bigr).
\end{eqnarray*}
This yields the result of the lemma, since the bias of the toy
estimator is given by
\[
\frac{\beta_1(u)b_n^2}{\sigma_1(u)}+o \bigl(b_n^2 \bigr), %
\]
which implies, by the preceding, that
\[
F_n^{\toy}(t)+d_{F_0}(t) \biggl\{\int
_{u=0}^t\frac{\gamma
_n(u)g(u,t)}{F_0(t)-F_0(u)} \,du+\int
_{u=t}^M\frac{\gamma
_n(u)g(t,u)}{F_0(u)-F_0(t)} \,du \biggr\} %
\]
satisfies (\ref{linearasympequation1}), apart from a term of order
$O_p(n^{-1/2})$.
\end{pf*}
\end{appendix}



%

\printaddresses
\end{document}